\documentclass[10pt, draft]{amsart}
\usepackage[english]{babel}
\usepackage{latexsym}

\usepackage[english]{babel}
\usepackage{latexsym}
\usepackage{amsmath}

\DeclareMathOperator{\ad}{ad}

\DeclareMathOperator{\spn}{span}

\DeclareMathOperator{\AFS}{\mathit{AFS}}
\DeclareMathOperator{\Der}{Der}
\DeclareMathOperator{\sL}{\mathit{s\ell}}

\newcommand{\F}{\mathbb F}
\newcommand{\E}{\mathbb E}

\newcommand{\N}{\mathbb N}
\newcommand{\Z}[0]{\mathbb Z}

\newcommand{\n}[0]{\text{$\underline{n}$}}

\newtheorem{dummy}{Dummy}

\numberwithin{dummy}{section}
\numberwithin{equation}{section}

\newtheorem{theorem}[dummy]{Theorem}

\theoremstyle{definition}
\newtheorem{definition}[dummy]{Definition}

\theoremstyle{remark}
\newtheorem{rem}[dummy]{Remark}
\newtheorem*{rem*}{Remark to ourselves}

\begin{document}

\bibliographystyle{amsalpha}
\author{Marina Avitabile}
\email{marina.avitabile@unimib.it}
\address{Dipartimento di Matematica e Applicazioni\\
  Universit\`a di Milano-Bicocca\\
  via Cozzi 53\\
  I-20125 Milano\\
  Italy}

\author{Sandro Mattarei}
\email{mattarei@science.unitn.it}
\address{Dipartimento di Matematica\\
  Universit\`a degli Studi di Trento\\
  via Sommarive 14\\
  I-38050 Povo (Trento)\\
  Italy}

\thanks{Both authors are members of INdAM-GNSAGA, Italy.
  They acknowledge financial  support from Ministero dell'Istruzione,
  dell'Universit\`a e della  Ricerca, Italy, to  the project
  ``Graded Lie algebras  and pro-$p$-groups of finite width''.}
\title[Thin loop algebras of {A}lbert-{Z}assenhaus algebras]{Thin loop algebras of {A}lbert-{Z}assenhaus algebras}
\date{\today}
\begin{abstract}
Thin Lie algebras are Lie algebras over a field, graded over the positive integers and
satisfying a certain narrowness condition.
In particular, all homogeneous components have dimension one or two, and are called {\em diamonds} in the latter case.
The first diamond is the component of degree one, and the second diamond can only occur in degrees
$3$, $5$, $q$ or $2q-1$, where $q$ is a power of the characteristic of the underlying field.
Here we consider several classes of thin Lie algebras  with second
diamond in degree $q$.
In particular, we identify the Lie algebras in one of these classes with suitable loop algebras of certain
Albert-Zassenhaus Lie algebras.
We also apply a deformation technique to recover other thin Lie algebras
previously produced as loop algebras of certain graded Hamiltonian Lie algebras.
\end{abstract}
\maketitle

\section{Introduction}\label{sec:introduction}

A graded Lie algebra
\begin{equation*}
L= \bigoplus_{i=1}^{\infty} L_i
\end{equation*}
over a field $\F$ is said to be {\em thin} (according to~\cite{CMNS})
if $\dim(L_1)=2$ and the following
\emph{covering property} holds:
\begin{equation}\label{eq:covering}
\text{$L_{i+1}=[u,L_1]$ for every $0 \neq u \in L_i$, for all $i\ge 1$.}
\end{equation}
The definition implies at once that $L$ is generated by $L_1$ as a Lie algebra,
and that every homogeneous component has dimension $1$ or $2$.
A homogeneous component of dimension $2$ is referred to as a \emph{diamond},
In particular $L_1$ is a diamond, the {\em first diamond}.
Another consequence of the definition is that the filtration associated with the grading
coincides with the lower central series of $L$.
Because of this coincidence the {\em degree} of a homogeneous element has sometimes been called
{\em weight} or {\em class} in the literature.
Our preference here goes to the first name.

The covering property has the following equivalent formulation:
every graded ideal $I$ of $L$ is located between two consecutive terms $L^i$ of the lower central series of $L$,
in the sense that $L^i\supseteq I\supseteq L^{i+1}$ for some $i$.
We refer to~\cite{CMNS} for the easy proof that the two definitions are equivalent.
Borrowing some terminology from the theory of pro-$p$-groups (see Chapter~12 of the book~\cite{L-GMcKay}),
which translates naturally to the present setting
(with {\em graded ideal} replacing {\em closed normal subgroup,} and so on)
we can describe thin Lie algebras as the (positively) graded Lie algebras of {\em width} two and {\em obliquity} zero.

Thin Lie algebras constitute a natural generalization of the {\em graded Lie algebras of maximal class}
studied in~\cite{CMN} and~\cite{CN}.
In fact, the latter coincide with the thin Lie algebras which have no diamonds except the first one.
It is convenient, however, to explicitly exclude Lie algebras of maximal class from the definition
of thin algebras.
Therefore, as a general rule we require a thin Lie algebra to have at least two diamonds.
A natural exception to this rule occurs in characteristic two, as we explain
in the final part of Section~\ref{sec:second_diamond}.
In this paper we restrict our
attention to infinite-dimensional thin Lie algebras.

It is known from~\cite{CMNS} and~\cite{AviJur} that the second diamond of an infinite-dimensional thin Lie algebra
of zero or odd characteristic
can only occur in degree $3$, $5$, $q$ or $2q-1$ where $q$ is a power of the characteristic $p$.
All these possibilities really do occur, and examples can be constructed by taking
(the positive parts of twisted) {\em loop algebras} of suitable finite-dimensional Lie algebras.
(We explain this terminology in Section~\ref{sec:Nottingham}.)
In particular, examples of thin Lie algebras with second diamond in degree $3$ or $5$ can be produced
as loop algebras of classical Lie algebras of type $A_1$ or $A_2$ (see~\cite{CMNS} and~\cite{Car:thin_addendum})
and arise also as the graded Lie algebras associated to the lower central series of
certain $p$-adic analytic pro-$p$ groups of the corresponding types (see~\cite{Mat:thin-groups}, and also~\cite{KL-GP}).

The remaining cases are typically modular, and all the examples constructed so far
involve finite-dimensional simple Lie algebras of Cartan type.
For simplicity, in the following discussion we tacitly assume the characteristic of the ground field to be odd.
However, in the various sections of the paper we will also discuss the changes occurring in the case of characteristic two.
Thin algebras with second diamond in degree $2q-1$ have been called {\em $(-1)$-algebras} in~\cite{CaMa:thin}.
It was shown there how to assign a {\em type} to each diamond of a $(-1)$-algebra.
The position of the diamonds, and their types, which
take value in the underlying field plus infinity, are a set of invariants for the $(-1)$-algebra
which determine it up to isomorphism.
A variety of cases arise, according as the diamond types are all infinite, all finite, or of both types,
and some of the corresponding algebras have been constructed in~\cite{CaMa:thin,CaMa:Hamiltonian,AviMat:-1},
as loop algebras of Hamiltonian algebras $H(2:\n;\Phi(\tau))^{(1)}$
(see Section~\ref{sec:Cartan} for a definition).
The case of $(-1)$-algebras with all diamonds of infinite type illustrates alone the abundance of examples,
as these algebras are in a bijective correspondence (described in~\cite{CaMa:thin}) with a certain (large) subclass of the
graded Lie algebras of maximal class.

An investigation of thin Lie algebras with the second diamond in degree $q$
was initiated in~\cite{CaMa:Nottingham}, where they were called {\em Nottingham} Lie algebras.
The reason for the name is that
the simplest example is the graded Lie algebra associated with the lower central series of the
Nottingham group (see \cite{Jenn,John,Cam}).
Here, too, one can assign a type to each diamond later than the first,
though by an essentially different rule from the case of $(-1)$-algebras.
As in that case, the types take value in the underlying field plus infinity,
but the second diamond of a Nottingham Lie algebra has always type $-1$.
Diamonds of type zero or one are really one-dimensional components,
but it is convenient to allow them in certain positions and dub them {\em fake}.
Nottingham Lie algebras in which the
distance between the second and the third diamond (if any) is greater than the distance
between the first and the second diamond were classified in \cite{Young:thesis},
for $p>3$.

It was shown in~\cite{CaMa:Nottingham} that if the third diamond
in a Nottingham Lie algebra (with $p>5$) occurs in degree $2q-1$
and is of finite type different from zero or one, then $L$ is uniquely determined
by a certain finite-dimensional graded quotient of it.
It turns out that the diamond types then follow an arithmetic progression, which is thus determined by the
type of the third diamond. In the two cases where the third diamond is fake,
a similar result holds once some extra conditions are prescribed.
We recall more details from~\cite{CaMa:Nottingham} in Section~\ref{sec:Nottingham}.
Existence of all such algebras was supported by computational evidence,
but~\cite{CaMa:Nottingham} contains only {\em uniqueness} results,
in the form of finite presentations for (central extensions of) them,
essentially encoding their structure up to the third diamond.

Some of the algebras described in~\cite{CaMa:Nottingham}, namely,
those where the arithmetic progression is the constant sequence $-1$,
were produced in~\cite{Car:Nottingham}
as (twisted) loop algebras of Zassenhaus algebras $W(1;n)$.
(That paper contains a uniqueness proof as well, and so does~\cite{Car:Zassenhaus-three}, which deals with the case of
characteristic three.)
Those where the type of the third diamond belongs to the prime field (and hence so do all diamond types)
but is different from  $-1$
were constructed in~\cite{Avi} as loop algebras of graded Hamiltonian algebras $H(2;(1,n))^{(2)}$
(suitably extended by an outer derivation in a couple of cases).
Note that the absence of fake diamonds in the former case, and the presence of diamonds of both types zero and one
in the latter, reflects in the dimension of the simple algebra used in the construction, namely
a power of $p$ and two less than a power of $p$, respectively (for $p>2$).

The main goal of the present paper is a construction for the Nottingham Lie algebras with third diamond in degree $2q-1$
and of finite type {\em not} belonging to the prime field.
The construction is given in Section~\ref{sec:finite} and is summarized in Theorem~\ref{thm:grading_finite}.
This step completes an existence proof for the whole range of thin Lie algebras considered in~\cite{CaMa:Nottingham},
namely, the Nottingham Lie algebras with third diamond of finite type, for $p>5$.
In this last case, again, there are no fake diamonds, and so a simple algebra of dimension a power of $p$ is involved.
In fact, here we construct such thin algebras as loop algebras of
Hamiltonian algebras of type $H(2;\n;\Phi(1))$.
These simple Lie algebras, together with the Zassenhaus algebras,
form the class of {\em Albert-Zassenhaus} algebras (see the end of Section~\ref{sec:Cartan}),
and this justifies the title of this paper.

We outline the contents of the various sections of the paper.
Sections~\ref{sec:Nottingham} and~\ref{sec:Cartan} are preliminary and essentially expository.
In Section~\ref{sec:Nottingham} we review some results from~\cite{CaMa:Nottingham}
and other papers on Nottingham Lie algebras
and give precise definitions for some concepts barely mentioned in this Introduction,
such as the diamond {\em types}.
Finally, we describe a {\em loop algebra} construction, which produces thin Lie algebras
starting from finite dimensional algebras with a certain cyclic grading.

Section~\ref{sec:second_diamond} contains a new approach to locating the second diamond of a thin Lie algebra.
In particular, we extend the definition of Nottingham Lie algebras
to ground fields of characteristic two.
The main peculiarity of this case is that the second diamond is itself fake, of type $-1\equiv 1$.
Thus, recognizing {\em the second diamond} as such in characteristic two
requires a revision of what we mean by this expression.
We do that in a uniform way for all thin Lie algebras which includes all characteristics.
All our results in Sections~\ref{sec:mixed}, \ref{sec:finite} and~\ref{sec:prime_field}
admit interpretation in the case of characteristic two, which we provide in separate remarks in each section for clarity.

We devote Section~\ref{sec:Cartan} to introducing certain simple modular Lie algebras
of Cartan type needed in this paper, namely, the Zassenhaus algebra $W(1;n)$ and the Hamiltonian algebras
$H(2;\n;\Phi)$ for various choices of the automorphism $\Phi$.
Most of the material in this section can be found in~\cite{Strade:book}, but we have
made an additional effort of including the case of characteristic two for Hamiltonian algebras.

Machine computations have suggested that there might be various Nottingham Lie algebras
where the third diamond has type $\infty$.
In fact, in Section~\ref{sec:mixed} we construct the simplest of those
as a loop algebra of $H(2;(1,n);\Phi(1))$ with respect to a certain grading.
This Nottingham Lie algebra has diamonds occurring at regular intervals,
the type pattern being isolated occurrences of diamonds of type $-1$ separated by sequences of
$r-1$ diamonds of type $\infty$, for some power $r$ of $p$.
The first author of this paper had previously given a construction for this algebra (unpublished),
but much more involved than the one given here.
Another possible diamond pattern emerging from machine computations has diamonds of finite types
separated by bunches of diamonds of infinite type, with the finite types following an arithmetic progression.
These algebras are under investigation by the first author of this paper.

In Section~\ref{sec:finite} we realize the main goal of the paper, by constructing
a Nottingham Lie algebra with any prescribed finite type of the third diamond outside the prime field.
These algebras have diamond occurring at regular intervals, with types following an arithmetic progression,
and are constructed as loop algebras of $H(2;(1,n);\Phi(1))$ with respect to a different cyclic grading
from that used in the previous section.

We have mentioned earlier that Nottingham Lie algebras with all diamond types in the prime field
have been constructed in other papers.
In Section~\ref{sec:prime_field} we obtain them all again starting from
those with diamond types finite but not all in the prime field produced in Section~\ref{sec:finite}.
Those of~\cite{Car:Nottingham} with all diamonds of type $-1$ are obtained from a Zassenhaus subalgebra of
$H(2;(1,n);\Phi(1))$ isomorphic with $W(1;n)$.
Dealing with the other case, originally studied in~\cite{Avi}, of non-constant finite types in the prime field,
involves viewing $H(2;(1,n);\Phi(1))$ as a deformation of its associated graded Lie algebra
with respect to the standard filtration,
which is a certain central extension of $H(2;(1,n))^{(1)}$.

\section{Nottingham Lie algebras}\label{sec:Nottingham}

Let $L=\bigoplus_{i=1}^{\infty}L_i$
be a thin Lie algebra over the field $\F$, as defined in the Introduction.
In this paper we will assume without further mention that $L$ is infinite-dimensional.
The possible location of the second diamond in a thin Lie algebra
was determined in~\cite{AviJur}, extending more specialized results in~\cite{CMNS},
and building on results of~\cite{CaMa:thin} and~\cite{CaJu:quotients}.
It was proved there that in a thin Lie algebra $L$ over a field of characteristic $p\neq 2$ the second diamond
can only occur in a degree $k$ of the form $3$, $5$, $q$, or $2q-1$,
where $q$ is a power of $p$.
We comment more on this in Section~\ref{sec:second_diamond}.
In particular, $3$ and $5$ are the only possibilities in characteristic zero
(which was already known from~\cite{CMNS}).
From now on we assume that $\F$ has positive characteristic $p$.
Note that the hypothesis $p\neq 2$ was not assumed in~\cite{AviJur}.
However, the proof depends on results in~\cite{CaJu:quotients} for $p\neq 2$, and in~\cite{Ju:quotients} for $p=2$.
Since the latter paper contains errors (only partially corrected in~\cite{JuYo:quotients}),
the argument in~\cite{AviJur} is inconclusive in characteristic two.
Apart from this problem, the case $p=2$ requires separate consideration in various arguments.
To mention just one instance, the quotient $L/L^k$ is metabelian
provided $p>2$ (as is proved in~\cite{CaJu:quotients}),
but need not be if $p=2$ (see~\cite{Ju:quotients,JuYo:quotients}).
Therefore, in the following discussion we assume that $\F$ has odd characteristic $p$, and postpone consideration
of the case of characteristic two to Section~\ref{sec:second_diamond}.

In this paper we restrict our attention to thin Lie algebras with second diamond in degree
a power $q$ of the characteristic.
We call such thin Lie algebras {\em Nottingham Lie algebras},
as we have anticipated in the Introduction.
We point out that in small characteristics there may be coincidences among the possibilities $3$, $5$, $q$ and $2q-1$
for the degree of the second diamond, which we discuss in detail in Section~\ref{sec:second_diamond}.
In particular, in characteristic five one may want to exclude from the class of Nottingham Lie algebras
the thin Lie algebra with second diamond in degree $5$ considered in~\cite{CMNS}, see Remark~\ref{rem:k=q=5}.
Although such a terminology choice is really a matter of convenience, that Lie algebra of~\cite{CMNS}
does not share several common features of Nottingham Lie algebras.
However, the situation gets much more complicated in the case of characteristic three and second diamond in degree $q=3$,
which comprises a host of thin Lie algebras, some of which we mention in Remark~\ref{rem:k=q=3}.
Because of this and more technical reasons in this section we assume $q>3$.

Thus, suppose $L$ is a thin Lie algebra with second diamond in degree $q>3$,
a power of the odd characteristic $p$.
Choose a nonzero element $Y\in L_1$ with $[L_2,Y]=0$.
According to~\cite{CaJu:quotients} we have
\begin{equation}\label{eq:first_chain}
C_{L_1}(L_2)=\cdots=C_{L_1}(L_{q-2})=\langle Y \rangle,
\end{equation}
where
$C_{L_1}(L_i)= \{ a \in L_1 \mid [a,b]=0 \textrm{ for every $b \in L_i$} \}$,
the centralizer of $L_i$ in $L_1$.
It follows from~\cite{Car:Nottingham} that one can choose
$X\in L_1\setminus\langle Y\rangle$ such that
\begin{equation}\label{eq:second_diamond_type}
[V,X,X]=0=[V,Y,Y], \quad [V,Y,X]=-2[V,X,Y],
\end{equation}
where $V$ is any nonzero element of $L_{q-1}$.
Here we adopt the left-normed convention for iterated Lie brackets, namely, $[a,b,c]$ denotes $[[a,b],c]$.
With respect to the notation of~\cite{Car:Nottingham,CaMa:Nottingham}
we have switched to capital letters for the generators $X$ and $Y$,
to avoid conflict of meaning with the
divided power indeterminates we use in later sections.

The proof of~\eqref{eq:second_diamond_type} given in~\cite{Car:Nottingham} is not very easy to locate,
being embedded in a rather complex proof that a certain Lie algebra is finitely presented.
Therefore, we isolate that part of the argument here for convenience.

\begin{proof}[Proof of~\eqref{eq:second_diamond_type}]
Here (but nowhere else in this paper) we need the identity
\[
[a,[b,\underbrace{c,\ldots,c}_{r}]]=\sum_{i=0}^r(-1)^i\binom{r}{i}
[a,\underbrace{c,\ldots,c}_{i},b,\underbrace{c,\ldots,c}_{r-i}],
\]
which holds for $a,b,c$ in an arbitrary Lie algebra and follows by induction from the Jacobi identity.
Because of~\eqref{eq:first_chain}, for an arbitrary choice of $X\in L_1\setminus\langle Y\rangle$
and setting $V=[Y,\underbrace{X,\ldots,X}_{q-2}]$ we have
\[
0=[[Y,\underbrace{X,\ldots,X}_{q-3}],[X,Y,Y]]=
[Y,\underbrace{X,\ldots,X}_{q-2},Y,Y]=[V,Y,Y],
\]
which is the second relation in~\eqref{eq:second_diamond_type}.
The covering property~\eqref{eq:covering} implies that $L_{q+1}$ has dimension one and is spanned by $[V,Y,X]$.
In particular, both $[V,X,X]$ and $[V,X,Y]$ are multiples of $[V,Y,X]$.
Similarly, we have
\begin{align*}
0&=(-1)^{(q-1)/2}[[Y,\underbrace{X,\ldots,X}_{(q-1)/2}],[Y,\underbrace{X,\ldots,X}_{(q-1)/2}]]=
\\
&=-\binom{(q-1)/2}{(q-3)/2}[V,Y,X]+[V,X,Y]=
\frac{1}{2}[V,Y,X]+[V,X,Y],
\end{align*}
which yields the third relation in~\eqref{eq:second_diamond_type}.
Hence $[V,X,Y]\neq 0$, and because
\begin{align}\label{eq:VXX=0}
[V,X+\alpha Y,X+\alpha Y]&=[V,X,X]+\alpha[V,X,Y]+\alpha[V,Y,X]
\\
\notag
&=[V,X,X]-\alpha[V,X,Y]
\end{align}
we can make the first relation in~\eqref{eq:second_diamond_type} hold by replacing
$X$ with $X+\alpha Y$ for a suitable scalar $\alpha$.
\end{proof}

It follows from~\eqref{eq:VXX=0}
that an element $X$ satisfying relations~\eqref{eq:second_diamond_type}
is determined up to a scalar multiple, as well as $Y$,
which is any nonzero element of $C_{L_1}(L_2)$.

Relations~\eqref{eq:second_diamond_type} also imply that $L_{q+1}$, the homogeneous component
immediately following the second diamond, has dimension (at most) one.
This extends to a more general fact about thin Lie algebras
(not necessarily with second diamond in degree $q$)
which is proved in~\cite{Mat:chain_lenghts}:
in a thin Lie algebra of arbitrary characteristic with $\dim(L_3)=1$,
two consecutive components cannot both be diamonds.
However, there are plenty of thin Lie algebras where both $L_3$ and $L_4$ are diamonds,
part of which are described in~\cite{GMY}, see Remark~\ref{rem:k=3}.

Suppose now that $L_{i}$ is a diamond of $L$ in degree $i>1$.
Let $V$ be a non-zero element in $L_{i-1}$, which has dimension one because of the previous paragraph,
whence $[V,X]$ and $[V,Y]$ span $L_i$.
If the relations
\begin{equation}\label{eq:diamond_type}
[V,Y,Y]=0=[V,X,X], \quad
(1- \mu)[V,X,Y]=\mu [V,Y,X],
\end{equation}
hold for some $\mu \in \F$, then we say that the diamond $L_i$ has (finite) type  $\mu$.
In particular, relations~\eqref{eq:second_diamond_type} postulate that the second diamond $L_q$ has type $-1$.
The third relation in~\eqref{eq:diamond_type} has a more natural appearance when written in terms of the generators
$Z=X+Y$ and $Y$, which were more convenient in~\cite{CaMa:Nottingham}, namely,
$[V,Z,Y]=\mu [V,Z,Z]$.
We also define $L_i$ to be a diamond of
type $\infty$ by interpreting the third relation in~\eqref{eq:diamond_type} as $-[V,X,Y]=[V,Y,X]$ in this case.
We stress that the type of a diamond of a Nottingham Lie algebra $L$ is
independent of the choice of the graded generators $X$ and $Y$.
In fact, we have seen above that our requirements $[L_2,Y]=0$ and~\eqref{eq:second_diamond_type}
determine $X$ and $Y$ up to scalar multiples,
and replacing them with scalar multiples does not affect~\eqref{eq:diamond_type}.
We note in passing that the situation is different for the thin Lie algebras with second diamond in degree $2q-1$,
where the diamond types are defined differently, and there is no canonical choice
for the type of the second diamond, see~\cite{CaMa:thin,CaMa:Hamiltonian,AviMat:-1}.

According to the defining equations~\eqref{eq:diamond_type}, a diamond of type $\mu=0$ should satisfy $[V,X,Y]=0=[V,X,X]$, while
a diamond of type $\mu=1$ should satisfy $[V,Y,X]=0=[V,Y,Y]$.
But then $[V,X]$, or $[V,Y]$, respectively, would be a central element in $L$ and,
therefore, would vanish because of the covering property~\eqref{eq:covering}.
Hence, $L_i$ would be one-dimensional after all, and hence not a genuine diamond.
Thus, strictly speaking, diamonds of type $0$ or $1$ cannot occur.
Nevertheless, we allow ourselves to call {\em fake diamonds} certain one-dimensional homogeneous components of $L$
and assign them type $0$ or $1$ if the corresponding relations hold, as described by~\eqref{eq:diamond_type}.
This should be regarded as a convenient piece of terminology rather than a formal definition:
whenever a one-dimensional component $L_i$ formally satisfies the relations of a fake diamond of type $1$,
the next component $L_{i+1}$ satisfies those of a fake diamond of type $0$.
However, as a rule, only one of them ought to be called a fake diamond,
usually because it fits into a sequence of (genuine) diamonds occurring at regular distances.

It was proved in~\cite{CaMa:Nottingham} that
\begin{equation}
\label{eq:second_chain}
 C_{L_1}(L_{q+1})= \cdots = C_{L_1}(L_{2q-3})=\langle Y \rangle,
\end{equation}
provided $p>5$.
In particular, $L$ cannot have a third diamond in degree lower than $2q-1$.
The assumption on the characteristic was weakened in~\cite[Proposition~5.1]{Young:thesis}
to $p>3$ and $q\neq 5$.
The failure of these assertions for $q=p=5$ is shown by
the thin algebra constructed in~\cite{CMNS} as a loop algebra of a simple Lie algebra of classical type $A_2$.
This has diamonds in all degrees congruent to $\pm 1$ modulo $6$ and, in particular,
has its third diamond in degree $7=2q-3$.
As noted in~\cite[Proposition~5.2]{Young:thesis},
the proof of~\cite[Theorem~1]{CMNS} for the case $p>5$
can be adjusted to prove that this is the only exception to~\eqref{eq:second_chain} for $q=5$.
In characteristic $p=3$ there are infinitely many counterexamples to~\eqref{eq:second_chain}.
All known ones have a third diamond in degree $2q-3$, and form a family $T(q,1)$
described without proof in~\cite[Section~6.3]{Young:thesis},
based on computational evidence.
In fact, Young described a larger family $T(q,r)$ (with both $q$ and $r$ powers of three) of exceptional thin Lie algebras
in characteristic three.
The paper~\cite{Mat:char_3} contains an explicit construction for $T(q,r)$
based on Hamiltonian algebras $H(2;(1,n))^{(2)}$, where $3^{n+1}=qr$.

The third genuine diamond does occur in degree higher than $2q-1$ for some of the Nottingham Lie algebras
identified in~\cite{CaMa:Nottingham} and then
constructed in~\cite{Avi} (which we consider again in Section~\ref{sec:prime_field}).
In fact, in two (families) of those Lie algebras the diamond in degree $2q-1$ is fake (of type zero or one),
and the third genuine diamond occurs in higher degree (see the next paragraph).
Young has classified in~\cite[Theorem~5.7]{Young:thesis} all Nottingham Lie algebras
with third genuine diamond in degree higher than $2q-1$, for $p>3$.
Besides the Lie algebras from\cite{CaMa:Nottingham,Avi} which we have just mentioned,
his list comprises three countable families
and an uncountable one.
Two of those countable families actually consist of soluble Lie algebras
{\em of finite coclass}, having just two diamonds.
An almost immediate consequence is a classification of the graded Lie algebras (generated by $L_1$)
of coclass two (for $p>3$):
any such Lie algebra either belongs to one of those two families, or is a central extension
of a graded Lie algebra of maximal class.
Finally, the uncountable family
is in a bijection with a subfamily of the graded Lie algebras of maximal class of~\cite{CMN,CN}.
This is quite remarkable in view of the fact,
proved in~\cite{CaMa:thin},
that thin Lie algebras
with second diamond in degree $2q-1$ also include an uncountable family
in a bijection with a subfamily of the graded Lie algebras of maximal class.
There appears to be no obvious connection between the two constructions.

Assuming that $L$ does have a genuine third diamond of finite type $\mu_3$ (hence different from $0$ and $1$)
in degree $2q-1$,
detailed structural information on $L$ was obtained in~\cite{CaMa:Nottingham}, for $p>5$.
(The case $\mu_3=-1$, which includes a graded Lie algebra associated with the Nottingham group,
had already been dealt with in~\cite{Car:Nottingham} for $p>3$
and~\cite{Car:Zassenhaus-three} for $p=3$.)
It was shown in~\cite{CaMa:Nottingham} that $L$ is uniquely determined, that the diamonds occur in each degree
congruent to $1$ modulo $q-1$, and that their types follow an arithmetic progression
(thus determined by $\mu_3$).
This was proved by showing that the homogeneous relations satisfied by $L$
up to degree $2q$ define a finite presentation for a central extension of $L$
(second central in case $\mu_3=-2$),
or even for $L$ itself if $\mu_3\not\in\F_p$.
Similar structural results on $L$ were also obtained in the exceptional cases $\mu_3=0,1$,
where $L_{2q-1}$ is fake (and also $L_{3q-2}$ in the former case).
However, here uniqueness of $L$ can only be established by assuming its structure
up to the fifth or fourth diamond, respectively
(that is, up to the third genuine diamond).
In fact, this is where the further Lie algebras described in Theorem~5.7 of~\cite{Young:thesis} enter the picture.

We have mentioned in the Introduction that existence results for some of the Nottingham Lie algebras
described in~\cite{CaMa:Nottingham} have been established in various papers.
Before providing some more details of these results we describe a construction which is crucial for them.

\begin{definition}\label{def:loop_algebra}
Let $S$ be a finite-dimensional Lie algebra over the field $\F$ with a cyclic grading $S=\bigoplus_{k\in\Z/N\Z} S_k$,
let $U$ be a subspace of $S_{\bar 1}$ and let $t$ be an indeterminate over $\F$.
The {\em loop algebra} of $S$ with respect to the given grading and the subspace $U$ is the Lie subalgebra of
$S\otimes\F[t]$ generated by $U\otimes t$.
\end{definition}

We omit mention of $U$ when it coincides with $S_{\bar 1}$, as will always be the case in this paper,
but an example where it does not occurs in~\cite{Mat:char_3}.
Note that elsewhere, for instance in~\cite{CaMa:Hamiltonian}, the loop algebra is defined as
$\bigoplus_{k>0}S_{\bar k}\otimes\F t^k$ (where $\bar k=k+N\Z$),
while with the present definition the loop algebra is only a subalgebra of that.
The two definitions coincide exactly when $U=S_{\bar 1}$ generates $S$ as a Lie algebra and
$S_{\bar 1}^{N+1}=S_{\bar 1}$.
This does occur in several cases of interest
but the two definitions are definitely not equivalent when $S$ is not perfect, for example.
Typically, $S$ has finite dimension, and the loop algebra construction
serves to produce an infinite-dimensional Lie algebra by ``replicating'' its structure periodically.
In the applications which we have in mind $U$ is two-dimensional
and $S$ is simple or close to being simple.
Proving that such a loop algebra is thin (in particular, the verification of the covering property~\eqref{eq:covering})
can conveniently be done inside the finite-dimensional Lie algebra $S$.
This is best illustrated by working out specific examples,
such as those in Sections~\ref{sec:mixed} and~\ref{sec:finite}.
Finally, note that if a Lie algebra $S$ is graded over an arbitrary cyclic group $G$ of order $N$ we need
to fix a specific isomorphism of $G$ with $\Z/N\Z$
(or, equivalently, a distinguished generator ``$\bar 1$'' for $G$)
for the corresponding loop algebra to be defined.

Nottingham Lie algebras with third diamond (in degree $2q-1$ and) of type $\mu_3=-1$
were constructed in~\cite{Car:Nottingham}, for every odd $p$
(but see also~\cite{Car:Zassenhaus-three} for the peculiarities of characteristic three),
as loop algebras of Zassenhaus algebras $W(1;n)$, which are simple Lie algebras of dimension a power $q$ of $p$,
with respect to a suitable grading.
Among these is the graded Lie algebra associated with the lower central series of the Nottingham group
(see the last section of the book~\cite{L-GMcKay}, for example), which occurs when $q=p$.

Nottingham Lie algebras with third diamond of type $\mu_3\in\F_p\setminus\{-1\}$
were constructed in~\cite{Avi}.
They have diamonds (some of which fake) in all degrees congruent to $1$ modulo $q-1$,
with types following an arithmetic progression.
Although the assumption $p>5$ stated in~\cite{Avi} is needed when dealing with finite presentations
(and to prove uniqueness results such as those of~\cite{CaMa:Nottingham} mentioned above),
the construction given there really works in every odd characteristic (see Remark~\ref{rem:k=q=3} for characteristic three)
and, with suitable adjustments, even in characteristic two (see Section~\ref{sec:second_diamond}).
For $\mu_3\neq -2,-3$ those algebras are loop algebras of graded Hamiltonian algebras $H(2;(1,n))^{(2)}$,
which are simple Lie algebras (for $p$ odd) of dimension two less than a power of $p$,
see Section~\ref{sec:Cartan}.
In the exceptional cases where $\mu_3=-2$ or $-3$, the very first diamond $L_1$ would be fake,
if the arithmetic progression of diamond types were extended backwards to include it.
Therefore, in order to obtain a two-dimensional generating subspace in degree one
it is necessary to extend $H(2;(1,n))^{(2)}$ by an outer derivation before applying the loop algebra construction.

In Section~\ref{sec:finite} we construct, over a perfect field $\F$ of arbitrary positive characteristic
(but see Section~\ref{sec:second_diamond} to make sense of characteristic two), Nottingham Lie algebras with
diamonds in all degrees congruent to $1$ modulo $q-1$,
with types following an arithmetic progression, and $\mu_3$ an arbitrary element of $\F$ not lying in the prime field.
They are loop algebras of Albert-Zassenhaus algebras $H(2;(1,n);\Phi(1))$,
which we define in Section~\ref{sec:Cartan}.
Finally, in Section~\ref{sec:prime_field} we apply a deformation argument to these Nottingham Lie algebras
to recover those of~\cite{Avi} considered in the previous paragraph, having diamond types in the prime field.

So far we have not considered the possibility that the third diamond in a Nottingham Lie algebra has infinite type.
We construct one such algebra in Section~\ref{sec:mixed}, again as a loop algebra of $H(2;(1,n);\Phi(1))$,
but with respect to a different grading.
There exist more Nottingham Lie algebras with third diamond of infinite type, and
they are being investigated in~\cite{AviJur:Nottingham_mixed}.

\section{The second diamond in a thin Lie algebra}\label{sec:second_diamond}

In characteristic two we have to adjust the definition of a Nottingham Lie algebra,
because the second diamond itself is fake, having type $-1=1$.
Thus, the second diamond is not recognised as such under the natural but restrictive definition of a diamond
as a two-dimensional component.
This causes the discrepancy, mentioned in the first paragraph of Section~\ref{sec:Nottingham},
between $L/L^k$ being metabelian for $p$ odd (which is proved in~\cite{CaJu:quotients}),
where $k$ is the degree of the second diamond,
and not being necessarily so for $p=2$ (see~\cite{Ju:quotients,JuYo:quotients}).
This discrepancy can be partially resolved by redefining the parameter $k$ which marks the degree of the second diamond
in a thin Lie algebra.
Before we do that we recall a concept from the theory of {\em graded Lie algebras of maximal class}
developed in~\cite{CMN,CN,Ju:maximal}.
To simplify the statements, without further mention we assume graded Lie algebra of maximal class
to be infinite-dimensional.
To each non-metabelian graded Lie algebra of maximal class $L$
one associates a certain parameter $q$, which is a power of $p$.
The usual way of defining this parameter is considering the smallest integer $k>2$ such that
$C_{L_1}(L_{k})\neq C_{L_1}(L_2)$.
According to Theorem~5.5 of~\cite{CMN}, $k$ has always the form $2q$ for some power $q$ of $p$,
and $q$ is called the {\em parameter} of $L$.
However, an equivalent way is setting $k:=\dim(L/L^{(2)})-1$, where $L^{(2)}=[[L,L],[L,L]]$
denotes the second derived subalgebra of $L$.
We take the latter approach to define a parameter $k$ for thin Lie algebras,
leaving aside the case where $\dim(L_3)=2$, which we consider separately in Remark~\ref{rem:k=3}.
It turns out that $k$ defined in this way coincides with the degree of the second diamond
in odd characteristic, but not always in characteristic two.
In the following theorem we essentially restate some known results in terms of this new definition of $k$.
For simplicity we restrict ourselves to the case of infinite-dimensional thin Lie algebras,
as everywhere in this paper, but we point out that the result remains true
for thin Lie algebras of finite dimension large enough (with respect to $k$).

\begin{theorem}\label{thm:parameter}
Let $L$ be an infinite-dimensional thin Lie algebra (hence not of maximal class, by convention) with $\dim(L_3)=1$, and set
$k=\dim(L/L^{(2)})-1$.
Then $k$ can only be $5$, $q$, or $2q-1$, where $q$ is some power of $p$.
Furthermore, $\dim(L_k)=1$ when $p=2$ and $k=q$,
and $\dim(L_k)=2$ otherwise.
\end{theorem}

\begin{proof}
Let $L_h$ be the second (genuine) diamond of $L$, hence $h>3$ by hypothesis.
Suppose first that $L/L_h$ is metabelian.
According to Theorem~2 of~\cite{AviJur} (which is not affected by the errors in~\cite{Ju:quotients} pointed out earlier),
$h$ can only be $5$, $q$, or $2q-1$.
Furthermore, an easy argument from~\cite[p.~283]{CMNS} shows that $h$ is odd, and that
$\dim(L^{(2)}\cap L_h)=1$.
Since $L$ is thin and $L^{(2)}$ is an ideal we have
$\dim(L^h/L^{(2)})=1$, and hence $k=h$.
We conclude that $k$ can only be $5$, $q$, or $2q-1$ in this case, with the exception of $q$ when $p=2$.

Now suppose that $L/L_h$ is not metabelian, hence $p=2$ according to~\cite{CaJu:quotients}.
We may deduce the conclusion from~\cite{Ju:quotients} as corrected in~\cite{JuYo:quotients},
but we prefer to proceed directly.
We have recalled above that Theorem~5.5 of~\cite{CMN}
determines the {\em length of the first constituent} in a graded Lie algebra of maximal class.
However, a closer look at its proof yields the following more precise result.
If the graded Lie algebra of maximal class $M$ has
exactly two components $M_i$ (besides $M_1$) which are not centralized by $C_{M_1}(M_2)$, say $M_k$ and $M_r$, with $k<r$,
then $k$ equals twice some power of the characteristic $p$.
In particular, when $p=2$ we conclude that $k$ is itself a power of two.
Since $L/L_h$ is of maximal class and not metabelian,
a suitable graded quotient $M$ of $L/L_{h+1}$ satisfies those assumptions, and the conclusion follows.
\end{proof}

In particular, we call a thin Lie algebra a {\em Nottingham Lie algebra}
if its parameter $k$ defined as in Theorem~\ref{thm:parameter} equals $q$, a power of $p$.
We will say that $L_k$ is the second diamond of $L$, regardless of the characteristic.
In particular, in characteristic two the second diamond is actually fake (that is, one-dimensional, of type $-1=1$),
and the second {\em genuine} diamond should perhaps be counted as the {\em third} diamond.
Apart from this peculiarity, the definition of diamond type given by relations~\eqref{eq:diamond_type}
and the related comments extend to the case of characteristic two.
These conventions will allow us to interpret Theorems~\ref{thm:grading_mixed} and~\ref{thm:grading_finite} for $p=2$.
We discuss further peculiarities of characteristic two in the final part of this section.

We mention for completeness that we call a thin Lie algebras a {\em $(-1)$-algebra}
(consistently with previous papers, and in absence of a better name)
if its parameter $k$ defined as in Theorem~\ref{thm:parameter} has the form $2q-1$, where $q$ is a power of $p$.

Generally speaking, as we have pointed out in the Introduction of the paper, when the characteristic $p$ is large enough
the possibilities $q$ and $2q-1$ for the degree of the second diamond correspond to thin Lie algebras
related with certain nonclassical finite-dimensional simple modular Lie algebras, as opposed to the values $3$ and $5$
of thin Lie algebras arising classical simple Lie algebras of type $A_1$ and $A_2$ as in~\cite{CMNS}.
When $p\le 5$, coincidences may arise between the former and latter cases.
In the following remarks we discuss those coincidences, as well as how to extend the definition of $k$ to values less than four.

\begin{rem}[$k=q=5$]\label{rem:k=q=5}
The first of the coincidences mentioned above occurs when $k=q=p=5$.
In view of the comments following Equation~\eqref{eq:second_chain},
it may be reasonable to extend the definition of Nottingham Lie algebras
to this case by adding the requirement that $\dim(L_7)=1$,
in order to exclude the loop algebras of the classical Lie algebras of type $A_2$ constructed in~\cite{CMNS}.
It is then possible to prove that~\eqref{eq:second_chain} holds in this case as well.
Nottingham Lie algebras fitting these prescriptions are those of~\cite{Car:Nottingham,Avi}
(see also Section~\ref{sec:prime_field}), as well as those of our Theorems~\ref{thm:grading_mixed} and~\ref{thm:grading_finite},
for $q=p=5$.
\end{rem}

\begin{rem}[$k=2q-1=5$]
A second diamond in degree $5$ also occurs for the $(-1)$-algebras with $q=p=3$.
In this case, however, the reason for the coincidence is deeper, as there is a thin Lie algebra fitting both recipes:
the classical Lie algebra of type $A_2$ is isomorphic with the Lie algebra
of Cartan type $H(2;(1,1),\Phi(\tau))^{(1)}$, the smallest {\em Block algebra} used in Theorem~7.2 of~\cite{CaMa:Hamiltonian}
(whose definition we recall in Section~\ref{sec:Cartan}).
Thus, in characteristic $3$ the Lie algebra with second diamond in degree $5$ constructed in~\cite{CMNS}
coincides with one of the infinite family of $(-1)$-algebras constructed in Theorem~7.2 of~\cite{CaMa:Hamiltonian}.
\end{rem}

\begin{rem}[$k=3$]\label{rem:k=3}
The case where the second diamond is $L_3$ was explicitly excluded from Theorem~\ref{thm:parameter}.
In fact, the approach taken there to defining the parameter $k$ cannot work in this case, because
the second derived subalgebra $L^{(2)}$ of a two-generated Lie algebra is always contained in the fifth term of the lower
central series.
A further obstacle here is that the characterization $C_{L_1}(L_2)=\langle Y\rangle$ of the distinguished
generator $Y$ fails, because $C_{L_1}(L_2)=\{0\}$.

This case comprises a host of thin Lie algebras with $\dim(L_4)=2$.
In particular, all metabelian thin Lie algebras have $\dim(L_3)=\dim(L_4)=2$,
and were shown in~\cite{GMY} to be in a bijective correspondence
with the quadratic extensions of the ground field $\F$.
Non-metabelian thin Lie algebras with all homogeneous components except $L_2$ having dimension two
can be produced via a quadratic field extension from certain graded Lie algebras of maximal class
of the type considered in~\cite{CV-L}.

Much more is known if we impose the requirement that $\dim(L_4)=1$.
In fact, if a thin Lie algebra in arbitrary characteristic satisfies $\dim(L_3)=2$ and $\dim(L_4)=1$,
commutation induces a nondegenerate symmetric bilinear form $L_1\otimes L_1\to L_4$ given by
$Z_1\otimes Z_2\mapsto [V,Z_1,Z_2]$, where $V$ is a nonzero element of $L_2$.
(Symmetry of the form is due to the fact that $L/L^5$ is metabelian for a two-generated Lie algebra,
and nondegeneracy follows from the covering property, see~\cite{CMNS}.)
Here we have two possibilities depending on the ground field,
according to the value of the discriminant of the associated quadratic form in $\F^\ast/(\F^\ast)^2$.
(This, besides those discussed in the previous paragraph,
is one of a few cases in the theory of thin Lie algebras where rationality questions arise.)
According to~\cite{CMNS}, for $p>3$ these two possibilities  lead to exactly two isomorphism classes
of (infinite-dimensional) thin Lie algebras, obtained as loop algebras of the split and the nonsplit form
of the classical Lie algebra of type $A_1$.
In the split case
there exist generators $X$ and $Y$ for $L$
such that $[V,X,X]=[V,Y,Y]=0$ and $[V,X,Y]=[V,Y,X]$, where $V=[Y,X]$.
\end{rem}

\begin{rem}[$k=q=3$]\label{rem:k=q=3}
In characteristic three the relations given above for the thin Lie algebra arising from $\sL_2$
(the split form of $A_1$) coincide with relations~\eqref{eq:second_diamond_type}
which describe the second diamond of a Nottingham Lie algebra.
There are infinitely many thin Lie algebras satisfying these requirements.
Besides the loop algebra of $\sL_2$ considered in~\cite{CMNS}, which is also
the graded Lie algebra associated with the Nottingham group (see~\cite{Car:Nottingham}) in view
of the isomorphism $\sL_2\cong W(1;1)$ in characteristic three,
we produce countably many in Theorem~\ref{thm:grading_mixed},
and a further one in Theorem~\ref{thm:grading_finite}.
In order that those theorems make sense for $q=p=3$,
we extend the definition of diamond type given by relations~\eqref{eq:diamond_type} to cover the present situation
(even though the characterization of $Y$ given by the condition $[L_2,Y]=0$ fails here).

Further thin Lie algebras of characteristic three with second diamond in degree three are the {\em Nottingham deflations} $N(3,r)$,
which have a diamond in degree one and in each degree congruent to $3$ modulo $3r-1$,
where $r$ is a power of $3$.
These are special cases of thin Lie algebras $N(q,r)$
described by Young in~\cite[Section~3.1]{Young:thesis} in arbitrary odd characteristic $p$,
for $q$ and $r$ powers of $p$, see also~\cite{Mat:char_3}.
Finally (as far as these authors know), there are algebras $T(3,r)$, belonging to a family $T(q,r)$
of exceptional thin Lie algebras of characteristic three (with $q$ and $r$ powers of $p=3$),
whose existence was conjectured in~\cite[Section~6.3]{Young:thesis}
on the basis of machine calculations and is proved in~\cite{Mat:char_3}.
The algebras $T(3,r)$
have a diamond in degree one and in each degree congruent to $3$ modulo $3r-3$.
\end{rem}

\begin{rem}[$k=2q-1=3$ and $k=q=2$]\label{rem:k=3or2}
The formulation of Theorem~\ref{thm:parameter} suggests that
when $\dim(L_3)=2$ and the characteristic is two we should allow the possibility
that $k=2$ and $L_2$ is a fake diamond.
In fact, certain thin Lie algebras in characteristic two with $\dim(L_3)=2$ admit
a double interpretation, both as $(-1)$-algebras, with second diamond in degree $2q-1=3$,
and as Nottingham Lie algebras, with (fake) second diamond in degree $q=2$.
We give an instance of this ambiguity in~Remark~\ref{rem:char_two_mixed_special}.
\end{rem}

We devote the rest of this section to a discussion of an interesting phenomenon in characteristic two.
The main peculiarity of characteristic two, which we mentioned earlier and made explicit in Theorem~\ref{thm:parameter},
is that the second diamond in degree $k=q$ is actually fake.
(However, its presence is revealed by the existence of a homogeneous element of $L_{[2]}$ in that degree;
here $L_{[p]}$ denotes the restricted closure of $L$ in $\Der(L)$.)
Consequently, our assumption to exclude graded Lie algebras of maximal class
from the definition of thin Lie algebras by requiring the existence of a second diamond becomes ambiguous
in characteristic two.
It seems to these authors that the notion of Nottingham (and hence thin) Lie algebras in characteristic two
should naturally embrace some of the graded Lie algebras of maximal class.
We do not intend to make here a definitive choice of terminology,
which is to some extent a matter of convenience,
but we do express some general comments in support of our view, supplemented by examples.

According to~\cite[Theorem~5.5]{CMN}, for an infinite-dimensional non-metabelian graded Lie algebra of maximal class $L$,
the earliest homogeneous component of $L$, after $L_1$,
which is not centralized by $\F Y=C_{L_1}(L_2)$, is $L_{2q'}$ for some power $q'$ of the characteristic.
In the terminology of~\cite{CN}, which revised that of~\cite{CMN} slightly,
the first constituent of $L$ has length $2q'$.
In characteristic two, a case dealt with in~\cite{Ju:maximal}, $q=2q'$ is itself a power of the characteristic.
If we choose $X$ and $V$ with $\F X=C_{L_1}(L_{q})$ and  $\F V=L_{q-1}$,
the relations~\eqref{eq:second_diamond_type} are satisfied, and
one is tempted to view $L$ as a thin Lie algebra with second diamond in degree $q$
(a Nottingham Lie algebra).
This becomes more compelling if $L$ has only two distinct {\em two-step centralizers}
(that is, each homogeneous component $L_i$ with $i\ge 2$ is centralized by either by $X$ or $Y$).
In fact, in that case one may place a fake diamond roughly in correspondence of each $L_i$ with $C_{L_1}\neq\F Y$,
namely, one may either call $L_i$ a diamond of type $1$, or $L_{i+1}$ a diamond of type $0$.
(We have discussed the ambiguity in defining diamonds of types zero and one
after introducing fake diamonds in Section~\ref{sec:Nottingham}.)
We stop this discussion at this point and present some specific examples.
We refer to Section~\ref{sec:Cartan} for the definitions of the Lie algebras of Cartan type involved
and the exceptional isomorphisms occurring in characteristic two.

The construction of~\cite{Car:Nottingham} based on the Zassenhaus algebra $W(1;n)$ can be interpreted in characteristic two
to produce a Nottingham Lie algebra with a diamond of type $-1\equiv 1$ in each degree $\equiv 1\pmod{q-1}$, where $q=2^n\ge 4$.
(We only point out that in characteristic two $W(1;n)$ is not simple, but its derived subalgebra $W(1;n)^{(1)}$ is;
this is one case where $S_{\bar 1}^{N+1}\neq S_{\bar 1}$ occurs in our Definition~\ref{def:loop_algebra} of loop algebra.)
Thus, all diamonds except the first are fake, and the thin Lie algebra is actually a graded Lie algebras of maximal class.
In the notation of~\cite{CN,Ju:maximal} it has sequence of constituent lengths $2q'(2q'-1)^\infty$, where $q'=q/2$,
and is therefore isomorphic with $AFS(n-1,n,n,2)$.
In fact, in characteristic two the Block algebra $H(2;(1,n-1);\Phi(\tau))^{(1)}$
used to produce $AFS(n-1,n,n,2)$ is isomorphic with
the simple Zassenhaus algebra $W(1;n)^{(1)}$.

Similarly, one can interpret in characteristic two also
the construction in~\cite{Avi}, which produces Nottingham Lie algebras with diamond types running over the prime field
as loop algebras of graded Hamiltonian algebras $H(2;(1,n))^{(2)}$.
In fact, in~\cite[Theorem~7]{Avi}
the diamonds occur at each degree congruent to $1$ modulo $q-1$, where $q=2^n$, and have types $0$ and $1$ alternately,
hence all fake.
Viewed as graded Lie algebra of maximal class, that algebra has
sequence of constituent lengths $(2q',2q'-2)^\infty$ with $q'=q/2$.
The theory of constituents in~\cite{CMN}
shows that the algebra is obtained by inflation from
another one with sequence of constituent lengths $(2q'',2q''-1)^\infty$
(obtained by halving the preceding ones), where $q''=q'/2=q/4$.
The latter is a loop algebra of $W(1;n-1)^{(1)}$, as in the previous paragraph.
This corresponds to the fact that in characteristic two
$H(2;(1,n))^{(2)}$ is not simple, but only differentiably simple and isomorphic with $W(1;n-1)^{(1)}\otimes\mathcal{O}(1;1)$.

\section{The Zassenhaus algebra and certain Hamiltonian algebras}\label{sec:Cartan}
We recall the definitions and the main properties of the Zassenhaus algebra and of certain Lie algebras of Cartan type
belonging to the Hamiltonian series.
We refer to Chapter~2 of the recent book of Strade~\cite{Strade:book} for a more thorough introduction.
A brief introduction to some of these algebras, aimed at the non-expert, was given in~\cite{CaMa:Hamiltonian}.
Here we depart from the notation used
in~\cite{CaMa:Hamiltonian}, which conformed to~\cite{BO,BKK,Kos:beginnings},
and convert to that of~\cite{Strade:book}, which has become standard.
However, we make some sign choices differently from~\cite{Strade:book}, conforming to~\cite{BO,BKK,Kos:beginnings}
instead, to preserve consistency with the papers~\cite{CaMa:Hamiltonian,Avi,AviMat:-1}.
We also make slight adjustments with respect to~\cite{Strade:book}
to accommodate the case of characteristic two, which is disregarded there.
The reader less experienced in the theory of modular Lie algebras may as well skip most of the following discussion
(but also lose much of the fun)
and view the Zassenhaus algebra $W(2;\n)$ and the Hamiltonian algebras
$H(2;\n)$, $H(2;\n;\Phi(\tau))$, $H(2;\n;\Phi(1))$
considered in this paper simply as the appropriate vector spaces with Lie brackets {\em defined} by
Equations~\eqref{eq:Zassenhaus}, \eqref{eq:Poisson_0}, \eqref{eq:Poisson_2},
and~\eqref{eq:Poisson_1}-\eqref{eq:Poisson_exception}, respectively.

Let $\F$ be an arbitrary field of prime characteristic $p$, and $m$ a positive integer.
(We will point out which statements need the assumption that $\F$ is algebraically closed.)
The algebra $\mathcal{O}(m)$ of divided powers in $m$ indeterminates (of unrestricted heights)
is the $\F$-vector space of formal $\F$-linear combinations of monomials
$x^{(a)}:=x_1^{(a_1)}\cdots x_m^{(a_m)}$,
with multiplication defined by
$x_i^{(a_i)} x_i^{(b_i)}=\binom{k+l}{k} x_i^{(a_i+b_i)}$,
and extended by linearity and by postulating commutativity and
associativity of the multiplication.
A natural $\Z$-grading $\mathcal{O}(m)=\bigoplus_{k\in\Z}\mathcal{O}(m)_k$ is given
by assigning to each monomial $x^{(a)}$ its total degree $|a|$.
Let $\mathcal{O}((m))$ be the completion of $\mathcal{O}(m)$ with respect to the topology
associated with the descending chain ideals defined by
$\mathcal{O}(m)_{(k)}=\bigoplus_{j\ge k}\mathcal{O}(m)_j$.

We partially order multi-indices $a,b\in\N^m$ by writing
$a\le b$ if and only if $a_i\le b_i$ for all $i$, and set
$|a|:=a_1+\cdots+a_m$.
Let $\n=(n_1,\ldots,n_m)$ be an $m$-tuple
of positive integers, and
set $\tau(\n):=(p^{n_1-1},\ldots,p^{n_m-1})$.
The monomials
$x^{(a)}$ with $0\le a\le\tau(\n)$
span a subalgebra of $\mathcal{O}(m)$ of dimension $p^{|n|}$,
denoted by $\mathcal{O}(m;\underline{n})$.
It is clearly isomorphic to the tensor product
$\mathcal{O}(1;n_1)\otimes\cdots\otimes\mathcal{O}(1;n_m)$
Actually, as an algebra $\mathcal{O}(m;\underline{n})$ is isomorphic with the tensor product of $|n|$
copies of the truncated polynomial ring $\F[x]/(x^p)$
(a copy for each $x_j^{(p^{k_j})}$, for $1\le j\le m$ and $0\le k_j<n_j$).
However, $\mathcal{O}(m;\underline{n})$ is endowed with an additional structure,
namely a set of {\em divided power maps}, which tie the various
$p$-(divided) powers of the same ``indeterminate'' together.
We refer to~\cite{Strade:book} for more details and only mention
that the definition of {\em special} derivations given in~\cite{Strade:book}
singles out exactly those (continuous) derivations of $\mathcal{O}(m)$ (resp.~$\mathcal{O}((m))$) which are compatible with
the divided power maps in a natural sense.
It turns out that the special derivations of $\mathcal{O}(m)$
(resp.~$\mathcal{O}((m))$, $\mathcal{O}(m;\underline{n})$) are
those of the form
$D=f_1\, \partial_1+\cdots+f_m\, \partial_m$
with $f_j\in\mathcal{O}(m)$ (resp.~$\mathcal{O}((m))$, $\mathcal{O}(m;\underline{n})$),
where $\partial_i(x_j^{(k)})=\delta_{ij}x_j^{(k-1)}$ (with $x_j^{(-1)}=0$), thus acting as
$Dx_j^{(k)}=f_jx_j^{(k-1)}$.
The special derivations of $\mathcal{O}(m)$ (resp.~$\mathcal{O}((m))$, $\mathcal{O}(m;\underline{n})$)
form a Lie algebra with respect to the Lie bracket, denoted by $W(m)$
(resp.~$W((m))$, $W(m;\underline{n})$) and called a {\em Witt algebra}.
(It is often referred to in the literature as the {\em generalized Jacobson-Witt algebra},
the name {\em the Witt algebra} being reserved to $W(1;1)$.)
It inherits a $\Z$-grading $W(m)=\bigoplus_{k\in\Z}W(m)_k$
(and a corresponding filtration) from the natural $\Z$-grading of $\mathcal{O}(m)$, where
$W(m)_{k}=\sum_{i=1}^m\mathcal{O}(m)_{k+1}\partial_i$.
In particular,
$W(m;\underline{n})=\sum_{i=1}^m\mathcal{O}(m;\underline{n})\partial_i$
has dimension $mp^n$.
It is a simple Lie algebra, except when $m=1$ and $p=2$ (see below).
Since in this paper we need at most two indeterminates, we find convenient to set $x:=x_1$ and $y:=x_2$.
We also use the standard shorthands $\bar{x}=x^{(p^{n_{1}}-1)}$ and  $\bar{y}=y^{(p^{n_{2}}-1)}$.
When $m=1$ we set $n:=n_1$, $\partial:=\partial_1$, and simply write $n$ for $\underline{n}$.

In this paper we only need the {\em Zassenhaus algebras}
$W(1;n)$, as they occur as distinguished subalgebras
of the Hamiltonian algebras which we will consider.
In this case the components of the standard grading are one-dimensional,
$W(1;n)_i$ being generated by $E_i= x^{(i+1)}\partial$, for $i=-1,\ldots,p^n-2$.
Direct computation shows that
\begin{equation}\label{eq:Zassenhaus}
[E_i,E_j]=
\left(
\binom{i+j+1}{j}-
\binom{i+j+1}{i}
\right)
E_{i+j}.
\end{equation}
In particular, we have $[E_{-1},E_j]=E_{j-1}$,
and $[E_0,E_j]=jE_j$.
The Zassenhaus algebra has also a grading over
(the additive group of) $\F_{p^n}$, with graded basis
consisting of the elements $e_{\alpha}$, for $\alpha\in\F_{p^n}$,
which satisfy
\begin{equation}\label{eq:Zassenhaus_group}
[e_{\alpha},e_{\beta}]=
(\beta-\alpha)
e_{\alpha+\beta}.
\end{equation}
One way of obtaining a basis $\{e_\alpha\}$ from the basis $\{E_i\}$
(see~\cite[Section~2]{CaMa:Hamiltonian} for details) is
\begin{equation*}
\begin{cases}
  e_0= E_{-1}+E_{p^n-2}
  &\\
  e_{\alpha}= E_{p^n-2}+ \sum_{i=-1}^{r} \alpha^{i+1}E_i & \text{for
    $\alpha\in\F_{p^n}^{\ast}$}.
\end{cases}
\end{equation*}
The former formula can be viewed as a special case of the
latter by setting $0^0:=1$.
In characteristic two the Zassenhaus algebra $W(1;n)$ is not simple, but its derived subalgebra
$W(1;n)^{(1)}=\langle E_i\mid i\not=p^n-2\rangle=
\langle e_{\alpha}\mid\alpha\not=0\rangle$
is simple.
The above transition formulas are clearly valid in this case as well.

The four series of Lie algebras {\em of Cartan type} $S$, $H$ and $K$ ({\em special,} {\em Hamiltonian} and {\em contact})
are defined as certain subalgebras of the generalized Jacobson-Witt algebra $W(m)$.
The {\em graded} representatives of each of these series, that is,
those which are graded subalgebras of $W(m)$ in its standard $\Z$-grading
(except in the contact case),
depend only on $m$, and $\n$ for the finite-dimensional ones.
However, the remaining ones depend also
on a choice of a certain continuous automorphism of $\mathcal{O}((m))$.
An alternative but equivalent description depends on a choice of a certain differential form.
The latter approach was taken in~\cite{CaMa:Hamiltonian}, following~\cite{BKK,Kos:beginnings} and other papers.
In this paper we consider only Hamiltonian algebras with $m=2$,
which can also be regarded as belonging to the series of Cartan type $S$.

As in~\cite[\S4.2]{Strade:book} we set
$H((2)):=
\{ D\in W((2))\mid
D(\omega_H)=0 \}$
and $H(2;\n):=H((2))\cap W(2;\n)$, where $\omega_H=dx\wedge dy$.
Here one may view the (formal) differential forms simply as elements of the
exterior algebra on the set $\{dx,dy\}$ over $\mathcal{O}((2))$.
In particular, the space of differential $2$-forms
is the free $\mathcal{O}((2))$-module
on the basis $\{dx\wedge dy\}$, and is a $W((2))$-module via
\[
D(f\,dx\wedge dy)=(Df)\,dx\wedge dy+f\,d(Dx)\wedge dy+f\,dx\wedge d(Dy)
\]
for $D\in W((2))$,
where $df=\partial_1(f)\,dx+\partial_2(f)\,dy$.
We define the linear map
\[
D_H:\mathcal{O}((2))\to W((2)),\quad
D_H(f):=\partial_2(f)\partial_1-\partial_1(f)\partial_2
\quad\text{for $f\in\mathcal{O}((2))$},
\]
which has kernel $\F 1$.
Here we have chosen this definition of $D_H$, which differs in sign from~\cite{Strade:book},
and agrees with with~\cite{Kos:beginnings} and~\cite{BO} instead,
for consistency with the papers~\cite{CaMa:Hamiltonian,Avi,AviMat:-1}.
It is then easy to see that
$H((2))=D_H\big(\mathcal{O}((2))\big)$ and, therefore,
\begin{align*}
H(2;\n)
&=D_H\big(\mathcal{O}(2;\n)\big)
\oplus\F x^{(p^{n_1}-1)}\partial_2
\oplus\F y^{(p^{n_2}-1)}\partial_1
\\
&=D_H\big(\mathcal{O}(2;\n)\oplus\F x^{(p^{n_1})}\oplus\F y^{(p^{n_2})}\big).
\end{align*}
Since the linear map $D_H$ is graded of degree $-2$,
$H(2;\n)$ is a graded subalgebra of $W(2;\n)$, that is,
$H(2;\n)=\bigoplus_{k\in\Z}H(2;\n)_k$, where
$H(2;\n)_k=H(2;\n)\cap W(2;\n)_k$.
Its second derived subalgebra
\[
H(2;\n)^{(2)}=\spn\{
D_H(x^{(i)}y^{(j)})\mid
0<(i,j)<\tau(\n)\}
\]
has dimension $p^{|n|}-2$, and is a simple Lie algebra if $p>2$.
(The Lie algebras
$H(2;\n)$ and $H(2;\n)^{(2)}$
were denoted by
$\tilde H(2:\mathbf{n};\omega_0)$ and
$H(2:\mathbf{n};\omega_0)$
in~\cite{CaMa:Hamiltonian}, where $\omega_0=\omega_H$.)
In characteristic two, $H(2;\n)^{(2)}$ is simple provided $n_1>1$ and $n_2>1$,
as one can prove along the lines of the proof of Theorem~3.5 or~Theorem~4.5 of~\cite[Chapter 4]{SF}.
However,
$H(2;(1,n))^{(2)}$ has
$\spn\{D_H(y^{(j)})\mid 0<j<2^{n_2}\}$
as an ideal.
In fact, $H(2;(1,n))^{(2)}$ is the split extension of a simple Zassenhaus algebra $W(1;n)^{(1)}$ by its adjoint module,
which is also isomorphic with $W(1;n)^{(1)}\otimes\mathcal{O}(1;1)$.
(This fact was relevant in the final part of Section~\ref{sec:second_diamond}.)
In particular, $H(2;(1,n))^{(2)}$ is differentiably simple, see~\cite{Bl:differentiably-simple}.

Since
\begin{equation*}
[D_H(f),D_H(g)]=
D_H(D_H(f)(g))
\quad\text{for $f,g\in\mathcal{O}((2))$},
\end{equation*}
the map $D_H$ is a homomorphism from the $H((2))$-module $\mathcal{O}((2))$
onto the adjoint module for $H((2))$, with kernel $\F 1$.
The associative (and commutative) algebra $\mathcal{O}((2))$
can be endowed with an additional structure of Lie algebra
with respect to the {\em Poisson bracket}
$\{f,g\}=D_H(f)(g)=
\partial_2(f)\partial_1(g)-\partial_1(f)\partial_2(g)$
(which is the opposite as in~\cite{Strade:book}),
and the map $D_H$ yields a Lie algebra isomorphism from
$\mathcal{O}((2))/\F 1$ onto $H((2))$.
It is often convenient to carry out computations in $\mathcal{O}((2))$ with the Poisson product rather than in $H((2))$.
In particular, we have
\begin{equation}\label{eq:Poisson_0}
\begin{split}
\{x^{(i)}
y^{(j)},
x^{(k)}
y^{(l)}\}
&=
x^{(i)}
y^{(j-1)}
x^{(k-1)}
y^{(l)}-
x^{(i-1)}
y^{(j)}
x^{(k)}
y^{(l-1)}
\\
&=
N(i,j,k,l)\,
x^{(i+k-1)}
y^{(j+l-1)},
\end{split}
\end{equation}
where
\[
N(i,j,k,l):=
\binom{i+k-1}{i}
\binom{j+l-1}{j-1}-
\binom{i+k-1}{i-1}
\binom{j+l-1}{j}.
\]
The formula shows that $\mathcal{O}((2))$ (and, hence, $H((2))\cong\mathcal{O}((2))/\F1$)
is $\Z^2$-graded by assigning degree $(i,j)$ to the monomial $x^{(i+1)}y^{(j+1)}$.
To save on notation we will view such monomials as elements of $H((2))$,
omitting the isomorphism $D_H$ and leaving implicit that we mean their cosets modulo $\F1$.
This $\Z^2$-grading of $H((2))$, which is induced by a natural $\Z^2$-grading of the larger Lie algebra $W((2))$,
will be for us more important than the standard $\Z$-grading, where $x^{(i+1)}y^{(j+1)}$ has degree $i+j$.

Following~\cite[\S6.3]{Strade:book} and referring to that for more details,
we define two more Hamiltonian subalgebras of $W(2;\n)$, denoted by
$H(2;\n;\Phi(\tau))^{(1)}$ and $H(2;\n;\Phi(1))$.
We actually need only the latter in the present paper, but the former occurs in~\cite{CaMa:Hamiltonian,AviMat:-1}
and we describe it here for completeness.
All Lie algebras of Cartan type $H$ (or, equivalently, $S$)
associated with the parameters $m=2$ and $\n$ can be obtained by intersecting
$\varphi\circ H((2))\circ\varphi^{-1}$ with $W(2;\n)$ (and then possibly taking subalgebras),
for suitable choices of a continuous automorphism
$\varphi$ of $\mathcal{O}((2))$.
According to~\cite[Theorem~6.3.10]{Strade:book}, over an algebraically closed field of characteristic
$p>3$ there are only four choices for $\varphi$ (three if $n_1=n_2$)
leading to non-isomorphic simple Hamiltonian algebras associated with $m=2$ and $\n$.
They are given by
$\mathrm{id}$ (the identity automorphism, which gives $H(2;\n;\mathrm{id})=H(2;\n)$, already considered earlier),
$\Phi(\tau)$, and $\Phi(l)$ for $l=1,2$, where
\begin{align*}
\Phi(\tau)(x)&:=x-x^{(p^{n_1})}y^{(p^{n_1}-1)},\quad
&\Phi(\tau)(y)&:=y,
\\
\Phi(1)(x)&:=x\exp(x^{(p^{n_1})}),\quad
&\Phi(1)(y)&:=y,
\end{align*}
where $\exp(x^{(p^{n_1})}):=\sum_{i=0}^{\infty}x^{(ip^{n_1})}$.
Again to ensure consistency with~\cite{CaMa:Hamiltonian}
we have defined $\Phi(\tau)$ slightly differently from~\cite{Strade:book}
(where $\Phi(\tau)(x)$ is set to equal $x+x^{(p^{n_1})}y^{(p^{n_1}-1)}$ instead),
but equivalently as long as the ground field $\F$ is algebraically closed.
The automorphism $\Phi(2)$ is obtained from $\Phi(1)$ by interchanging the roles of $x$ and $y$,
and thus becomes unnecessary if we relax the convention of~\cite[Theorem~6.3.10]{Strade:book} that $n_1\le n_2$.

We deduce from~\cite[\S6.3]{Strade:book}, in particular from Equation~(6.3.4), that
\begin{align*}
H(2;\n;\Phi(\tau))
&=\{(1+\bar x\bar y)D\mid D\in H(2;\n)\}
\\
&=H(2;\n;\Phi(\tau))^{(1)}
\oplus\F y^{(p^{n_2}-1)}\partial_1
\oplus\F x^{(p^{n_1}-1)}\partial_2.
\end{align*}
According to~\cite[Theorem~6.3.10]{Strade:book}, its derived subalgebra
\[
H(2;\n;\Phi(\tau))^{(1)}=
\F(1+\bar x\bar y)\partial_1\oplus\F(1+\bar x\bar y)\partial_2\oplus
\bigoplus_{k\ge 0}{(H(2;\n)^{(1)})}_k
\]
has dimension $p^{|n|}-1$, and is a simple Lie algebra in every positive characteristic.
(For the case of characteristic $2$, which is excluded in~\cite[Theorem~6.3.10]{Strade:book},
we may refer to~\cite{CaMa:Hamiltonian}, where
$H(2;\n;\Phi(\tau))$ and $H(2;\n;\Phi(\tau))^{(1)}$
were denoted by
$\tilde H(2:\mathbf{n};\omega_2)$ and
$H(2:\mathbf{n};\omega_2)$.)
In characteristic two, $H(2;(1,n);\Phi(\tau))^{(1)}$
is isomorphic with the Zassenhaus algebra $W(1;n+1)^{(1)}$,
an isomorphism being obtained by mapping
$(1+\bar x\bar y)D_H(xy^{(j)})\mapsto E_{j-1}$ and $(1+\bar x\bar y)D_H(y^{(j)})\mapsto E_{j+2^{n}-2}$.

Similarly, according to~\cite[Theorem~6.3.10]{Strade:book},
\[
H(2;\n;\Phi(1))
=\spn\{
D_H(x^{(i)}y^{(j)})-\bar xx^{(i)}y^{(j)}\partial_2\mid
0\le (i,j)\le\tau(\n)\}
\]
has dimension $p^{|n|}$, and is a simple Lie algebra for $p>2$.
In characteristic two $H(2;\n;\Phi(1))$ is not simple.
In fact, its derived subalgebra is given by the same description given above for $H(2;\n;\Phi(1))$
but with $0\le (i,j)<\tau(\n)$.
It is easy to prove that $H(2;\n;\Phi(1))^{(1)}$ is simple, of dimension $2^{|n|}-1$.
(The proof of~\cite[Theorem~6.5.8]{Strade:book}, where the assumption $p>2$ is only needed
when $m=2r>2$, can be interpreted to cover this case.)

It is convenient to realize
$H(2;\n;\Phi(\tau))^{(1)}$ and $H(2;\n;\Phi(1))$
as $\mathcal{O}(2;\n)/\F 1$ or $\mathcal{O}(2;\n)$
with respect to appropriate Lie brackets.
We follow~\cite[\S6.5]{Strade:book}, taking into account the sign changes adopted here.
As in the case of $H(2;\n)$ one defines a certain linear map
$D_{H,\omega}:\mathcal{O}((2))\to W((2))$ with kernel $\F 1$,
which now depends on a more general differential form $\omega=\varphi(\omega_H)$,
where $\varphi$ is the continuous automorphism of $\mathcal{O}((2))$ associated to the Hamiltonian algebra
under consideration.
The forms $\omega$ associated with $H(2;\n)$, $H(2;\n;\Phi(\tau))$ and $H(2;\n;\Phi(1))$
are thus
\[
\omega_H=dx\wedge dy,\quad
(1-\bar x\bar y)dx\wedge dy,\quad
d\big(\exp(x^{p^{n_1}})x\,dy\big)=\exp(x^{p^{n_1}})dx\wedge dy,
\]
respectively.
(These differ by a coefficient $2$ from those in~\cite[\S6.5]{Strade:book}, which is harmless in odd characteristic,
but makes our formulas work also in characteristic two.
Another possible choice for the third form, made in~\cite{Kos:beginnings,BKK,Skr:Hamiltonian},
which is equivalent but less convenient for our computations, is $d\big(\exp(x)x\,dy\big)=(1+x)\exp(x)dx\wedge dy$.)
It turns out that $D_{H,\omega}$ maps $\mathcal{O}((2))$ onto
\[
H((2;\omega)):=
\{ D\in W((2))\mid
D(\omega)=0\}=
\varphi\circ H((2))\circ\varphi^{-1}
\]
and $H(2;\n;\varphi)$ is obtained as $H((2;\omega))\cap W(2;\n)$.
Since
\begin{equation*}
[D_{H,\omega}(f),D_{H,\omega}(g)]=D_{H,\omega}(D_{H,\omega}(f)(g))
\quad\text{for $f,g\in\mathcal{O}((2))$}
\end{equation*}
according to~\cite[Equation~(6.5.5)]{Strade:book},
the map $D_{H,\omega}$ becomes a Lie algebra homomorphism if we endow $\mathcal{O}((2))$
with the {\em Poisson bracket}
$\{f,g\}=D_{H,\omega}(f)(g)$.

When $\varphi=\Phi(\tau)$ we have
$D_{H,\omega}(f)=(1+\bar x\bar y)D_H(f)$ for all $f\in\mathcal{O}(2)$
(see~\cite[Equation~(6.5.4)]{Strade:book}).
It follows that
$\{f,g\}=(1+\bar x\bar y)\big(\partial_2(f)\partial_1(g)-\partial_1(f)\partial_2(g)\big)$.
In particular, $H(2;\n,\Phi(\tau))$ is isomorphic with the quotient of
$\mathcal{O}(2;\n)\oplus\F x^{(p^{n_1})}\oplus\F y^{(p^{n_2})}$
modulo its centre $\F 1$, with the Poisson bracket
\begin{equation}\label{eq:Poisson_2}
\{x^{(i)}
y^{(j)},
x^{(k)}
y^{(l)}\}
=
N(i,j,k,l)\,
(1+\bar x\bar y)x^{(i+k-1)}
y^{(j+l-1)}.
\end{equation}
The formula shows that $H(2;\n,\Phi(\tau))$ (despite not being $\Z^2$-graded, unlike $H(2;\n)$)
has a natural grading over the group $\Z^2/\Z\tau(\n)$, where the monomial $x^{(i+1)}y^{(j+1)}$
has degree $(i,j)+\Z\tau(\n)$.
This fact was exploited in~\cite{CaMa:Hamiltonian} to produce various cyclic gradings of $H(2;\n,\Phi(\tau))$
which give rise to graded Lie algebras of maximal class and thin Lie algebras by a loop algebra procedure.

When $\varphi=\Phi(1)$ we have
$D_{H,\omega}(ef)=D_H(f)-\bar x f\partial_2$ for all $f\in\mathcal{O}((2))$, where we have set $e:=\exp(x^{(p^{n_1})})$,
see~\cite[Equation~(6.5.8)]{Strade:book}.
Consequently, the map $D_{H,\omega}$ maps $e\mathcal{O}(2)$ onto $H(2;\omega):=H((2;\omega))\cap W(2)$.
In particular, $D_{H,\omega}$ maps $e\mathcal{O}(2;\n)$ onto $H(2;\n,\Phi(1))$.
It also follows that
\[
D_{H,\omega}(ef)(eg)=e\big(D_H(f)(g)+\bar x(\partial_2(f)g-f\partial_2(g)\big).
\]
Here we choose to define a Lie bracket by setting
$\{f,g\}:=e^{-1}D_{H,\omega}(ef)(eg)$
(which is different from the Poisson bracket defined above on $\mathcal{O}((2))$ with respect to a generic form $\omega$,
despite the same notation).
Then $H(2;\n,\Phi(1))$ becomes isomorphic to
$\mathcal{O}(2;\n)$ with respect to this Lie bracket.
Restricted to $\mathcal{O}(2;\n)$, that is, for $(0,0)\le (i,j),(k,l)\le\tau(\n)$,
the Lie bracket can be written as
\begin{align}
&\{x^{(i)}
y^{(j)},
x^{(k)}
y^{(l)}\}
=
N(i,j,k,l)\,
x^{(i+k-1)}
y^{(j+l-1)}
\quad\text{if $i+k>0$, and}\label{eq:Poisson_1}
\\
&\{y^{(j)},
y^{(l)}\}
=
\left(
\binom{j+l-1}{l}-
\binom{j+l-1}{j}
\right)\bar x
y^{(j+l-1)}.\label{eq:Poisson_exception}
\end{align}
It follows from Equations~\eqref{eq:Poisson_1} and~\eqref{eq:Poisson_exception} that the algebra $H(2;\n;\Phi(1))$
is graded over the group $\Z/p^{n_1}\Z\times\Z$
by assigning degree $(i+p^{n_1}\Z,j)$ to the monomial $x^{(i+1)}y^{(j+1)}$.

\begin{rem}
There is a way of making the special case~\eqref{eq:Poisson_1}
more similar to the general case~\eqref{eq:Poisson_exception}
by replacing $N(i,j,k,l)$ with
$
N'(i,j,k,l):=
\binom{i+k-1}{i}
\binom{j+l-1}{l}-
\binom{i+k-1}{k}
\binom{j+l-1}{j}.
$
In fact, $N'(i,j,k,l)=N(i,j,k,l)$ for $i,j,k,l\ge 0$ except when either $i=k=0$ or $j=l=0$,
which is immaterial in the former case and yields the nice formula
$\{y^{(j)},
y^{(l)}\}
=
N'(0,j,0,l)\bar x
y^{(j+l-1)}$
in the latter case.
Note also that replacing $N(i,j,k,l)$ with $N'(i,j,k,l)$ would leave
Equations~\eqref{eq:Poisson_0} and~\eqref{eq:Poisson_2} unaffected.
\end{rem}

The Hamiltonian algebras considered in this section are isomorphic to Lie algebras
produced by Albert and Frank in~\cite{AF} by means of different constructions
(see also~\cite[(V.4)]{Sel}).
In particular, an Albert-Zassenhaus algebra is a Lie algebra with a basis $\{u_\alpha\mid\alpha\in G\}$,
where $G$ is an additive subgroup of $\F$, and multiplication given by
\[
[u_\alpha,u_\beta]=(\beta-\alpha+\alpha\Theta(\beta)-\beta\Theta(\alpha))u_{\alpha+\beta}
\]
for some homomorphism $\Theta$ of $G$ into the additive group of $\F$.
This construction generalizes the construction~\eqref{eq:Zassenhaus_group} for the Zassenhaus algebra.
According to~\cite[Corollary~5.1]{BIO:Albert-Zassenhaus}, the Albert-Zassenhaus algebras of dimension $p^n$
over an algebraically closed field of characteristic $p>3$ are precisely
the Zassenhaus algebras $W(1;n)$ and the Hamiltonian algebras $H(2;\n;\Phi(1))$.
To justify the title of this paper, our main goal is
to construct certain cyclic gradings and corresponding loop algebras of the latter.

We stress that over an arbitrary field $\F$ of prime characteristic we take
$W(m;\n)$, $H(2;\n)$, $H(2;\n;\Phi(\tau))$ and $H(2;\n;\Phi(1))$
to denote the specific forms of these algebras as defined above,
although they may admit other forms when $\F$ is not algebraically closed.

\section{Nottingham Lie algebras with diamonds of types $-1$ and $\infty$}\label{sec:mixed}

The natural $(\Z/p^{n_1}\Z\times\Z)$-grading of $H(2;\n;\Phi(1))$ considered in the previous section
gives rise to a cyclic grading by viewing the degrees modulo
the subgroup of $\Z/p^{n_1}\Z\times\Z$ generated by $(-1,p^{n_2}-1)$.
More precisely, we obtain a grading $L=\bigoplus_{k\in\Z/N\Z}L_k$ of $L=H(2;\n;\Phi(1))$,
where $N=p^{n_1}(p^{n_2}-1)$,
by assigning degree $-(p^{n_2}-1)i-j+N\Z$ to the monomial $x^{(i+1)}y^{(j+1)}$.
In particular, $L_{\bar 1}$ is spanned by $x$ and $\bar y$, and all components of the grading are one-dimensional
except those with $k\equiv 1\pmod{p^{n_2}-1}$, which are two-dimensional.
This grading is the exact analogue of the grading of $H(2;\n;\Phi(\tau))^{(1)}$
described in~\cite[Subsection~5.2]{CaMa:Hamiltonian}.
We claim that the corresponding loop algebra is thin.
This follows from the following computations:
\begin{align*}
  & \{x^{(i)} y^{(j)}, x\} = x^{(i)}
  y^{(j-1)},\\
  & \{y^{(j)}, \bar y\} =
\begin{cases}
    -\bar x y^{(p^{n_2}-2)} &\text{if $j=0$,}
    \\
    2\bar x\bar y &\text{if $j=1$,}
    \\
    0 &\text{otherwise,}
  \end{cases}
\intertext{and, for $i>0$,}
  & \{x^{(i)} y^{(j)}, \bar y\} =
  \begin{cases}
    -x^{(i-1)} y^{(p^{n_2}-2)} &\text{if $j=0$,}
    \\
    x^{(i-1)} \bar y &\text{if $j=1$,}
    \\
    0 &\text{otherwise.}
  \end{cases}
\end{align*}
In fact, if $p>2$, which we assume for now,
the computations easily show that $L_{\bar 1}$ generates $H(2;\n;\Phi(1))$, and that
\begin{equation}\label{eq:covering_finite}
\text{$L_{k+\bar 1}=\{u,L_{\bar 1}\}$ for every $0 \neq u \in L_k$, for all $k\in\Z/N\Z$.}
\end{equation}
(We give more details about this in the next paragraph.)
This implies at the same time that the loop algebra of $L$ is thin,
and that it coincides with $\bigoplus_{k>0}L_{\bar k}\otimes\F t^k$.

Condition~\eqref{eq:covering_finite} expresses an analogue in $L$ of the covering property~\eqref{eq:covering}
for its loop algebra.
We now
check the validity of~\eqref{eq:covering_finite} in the only nontrivial cases,
namely, when $L_k$ is two-dimensional.
With a harmless abuse of notation with respect to Section~\ref{sec:Nottingham},
we denote by $X$ and $Y$ the elements $x$ and $\bar y$, respectively,
rather than the corresponding elements $x\otimes t$ and $\bar y\otimes t$ in the loop algebra.
In particular, the defining condition $\{L_{\bar 2},Y\}=0$ for $Y$ is satisfied.
Setting $V=x^{(i)}y$ for $0\le i<p^{n_1}$
we see that $L_{\bar k}$ for $k=-(p^{n_2}-1)(i-1)+1$
is spanned by $\{V,X\}=x^{(i)}$ and
\[
\{V,Y\}=
\begin{cases}
x^{(i-1)}\bar y&\text{if $i>0$,}\\
2\bar x\bar y&\text{if $i=0$.}
\end{cases}
\]
It follows that $\{V,X,X\}=\{V,Y,Y\}=0$ for all $i$.
These two relations together with the nonvanishing of both $\{V,X,Y\}$ and $\{V,Y,X\}$
imply the validity of the covering property for $L_{\bar k}$.
Therefore, the loop algebra of $L$ is thin.
Furthermore, since its second diamond occurs in degree $q$, it is a Nottingham Lie algebra,
and the definition of diamond type given by~\eqref{eq:diamond_type} applies.
Because $-\{VXY\}=\{V,Y,X\}$ for $i>0$ and
$2\{VXY\}=-\{V,Y,X\}$ for $i=0$,
the corresponding diamond has type $\infty$ in the former case and $-1$ in the latter.
For the purpose of reference we state our result formally.
\begin{theorem}\label{thm:grading_mixed}
Let $\F$ have odd characteristic $p$, and let $q=p^{n_2}$ and $r=p^{n_1}$, where $n_1,n_2>0$.
A grading of $H(2;\n;\Phi(1))$ over the integers modulo $(q-1)r$ is defined by declaring
\[
X=x
\qquad\textrm{and}\qquad
Y=\bar y
\]
to have degree one.
These elements generate $H(2;\n;\Phi(1))$ and are part of a graded basis
given by the monomials $x^{(i)}y^{(j)}$, where
$x^{(i)}y^{(j)}$ has degree
$(1-q)i-j+q$ modulo $(q-1)r$.
The corresponding loop algebra $N$ is thin with second diamond in degree $q$.
The diamonds occur in all degrees congruent to $1$ modulo $q-1$.
They have all type $\infty$, except those in degrees congruent to $q$ modulo $(q-1)r$, which have type $-1$.
\end{theorem}

\begin{rem}\label{rem:char_two_mixed}
A similar result as Theorem~\ref{thm:grading_mixed} holds in characteristic two,
but then $X=x$ and $Y=\bar y$ generate the derived subalgebra
$H(2;\n;\Phi(1))^{(1)}$, and the diamonds of type $-1$ in the loop algebra $N$ are fake,
as we have explained in Section~\ref{sec:Nottingham}.
The proof given above goes through with $L=H(2;\n;\Phi(1))^{(1)}$.
This Nottingham algebra $N$ falls into one of the cases described in~\cite{JuYo:quotients},
where the largest quotient of maximal class of a thin Lie algebra over a term of the lower central series
is not metabelian, but has constituent sequence $2q',2q'-2$, where $q'=q/2$.
Its existence provides a negative answer to the former of the two open questions
posed in the Introduction of~\cite{JuYo:quotients}
(but see also Remark~\ref{rem:char_two_finite}).
\end{rem}

\begin{rem}\label{rem:char_two_mixed_special}
The special case where $p=q=r=2$ is quite interesting.
The Hamiltonian algebra $H(2;(1,1);\Phi(1))^{(1)}$
is then the three-dimensional simple Lie algebra, and its loop algebra $N$
has a genuine diamond in each odd degree (and a fake one in each even degree).
Because of the many realizations of the three-dimensional simple Lie algebra,
this particular thin Lie algebra $N$ has several other interpretations besides that
as a Nottingham algebra with second diamond in degree $q=2$ which we have just given,
as we have anticipated in Remark~\ref{rem:k=3or2}.
For example, it is the loop algebra of $H(2;(1,1);\Phi(\tau))^{(1)}$
described in~\cite[Subsection~5.2]{CaMa:Hamiltonian}
(hence a maximal subalgebra of the graded Lie algebra of maximal class $\AFS(1,2,2,2)$, see~\cite{CMN,CN});
in this interpretation, its second diamond is (genuine) in degree $2q-1=3$,
and all its diamond have {\em infinite} type.
(Diamond types for $(-1)$-algebras are defined in a different way than for Nottingham Lie algebras,
see~\cite{CaMa:Hamiltonian}, for example.)
However, diamond types are ambiguous in this situation, and $N$ can also be viewed
as a loop algebra of $H(2;(1,1);\Phi(\tau))^{(1)}$
with respect to the different grading given in~\cite[Theorem~7.2]{CaMa:Hamiltonian},
thus, with all diamonds having {\em finite} types.
Finally, $N$ can be viewed as the case of characteristic two of the construction given in~\cite[Section~4]{CMNS},
namely, of a thin Lie algebra with second diamond in degree $3$ (independently of the characteristic).
In fact, the  three-dimensional Lie algebra $\mathcal{T}$ considered there is simple
in all characteristics but is classical of type $A_1$ only for $p$ odd
(unavoidably so, since the classical simple Lie algebra of type $A_1$
becomes nilpotent in characteristic two).
\end{rem}

\section{Nottingham Lie algebras with diamonds of finite types}\label{sec:finite}
In this section we produce a Nottingham Lie algebra with diamonds of finite types and not all in the prime field,
as the loop algebra of $H(2;(1,n);\Phi(1))$ with respect to a certain cyclic grading.
We follow a procedure used in~\cite[Section~7]{CaMa:Hamiltonian}
to construct $(-1)$-algebras with diamonds of finite types
as loop algebras of $H(2;(1,n);\Phi(\tau))^{(1)}$.
Differently from the construction in Section~\ref{sec:mixed}, here the ground field $\F$
is clearly not arbitrary, as it is bound to contain at least the various diamond types.
It will turn out after Theorem~\ref{thm:grading_finite} that the present construction works
over any perfect field $\F$ of characteristic $p>0$ which contains the diamond types.
However, it will be easier not to set that in advance, and simply allow ourselves to enlarge $\F$
when necessary (or just assume $\F$ algebraically closed for now).
As in Section~\ref{sec:mixed}, we assume that $p$ is odd,
and describe the necessary adjustments in characteristic two after Theorem~\ref{thm:grading_finite}.

We first work with
$L=H(2;\n;\Phi(1))$ and postpone the restriction $n_1=1$ to a later stage.
Consider the grading
$L=\bigoplus_{k\in\Z}\bar L_k$ of $L$ given by setting
\begin{equation*}
\bar L_{1-j}=
\langle
x^{(i)}y^{(j)}
\mid
i=0,\ldots,p^{n_1}-1
\rangle
\end{equation*}
for $j=0,\ldots,p^{n_2}-1$.
We identify the component $\bar L_0$ with the Zassenhaus algebra $W(1;n_1)$ via $E_i=x^{(i+1)}y$, for
$i=-1,\ldots,p^{n_1}-2$.
This is justified by
\begin{equation*}
\{x^{(i)}y,
x^{(k)}y\}
=
\left(
\binom{i+k-1}{k-1}-
\binom{i+k-1}{i-1}
\right)
x^{(i+k-1)}y.
\end{equation*}
We have mentioned in Section~\ref{sec:Cartan} that $W(1;n_1)$ has a grading
over the additive group of the finite field $\F_{p^{n_1}}$,
see Equation~\eqref{eq:Zassenhaus_group}.
More generally, we can obtain a grading
over a certain additive subgroup of
$\F$ (assumed large enough) as the root space decomposition of $W(1;n_1)$
with respect to an appropriate one-dimensional Cartan subalgebra, say $\langle e_0\rangle$,
where $e_0= y+\pi\bar x y$ for some nonzero coefficient $\pi$.
(The choice of $\pi$ was immaterial in~\cite{CaMa:Hamiltonian}, where we had $\pi=1$.)
Let $\sigma\in\F$ such that $\sigma^{p^{n_1}-1}=\pi$.
Then, nonzero elements in the corresponding one-dimensional root spaces are given by the formulas
\begin{equation*}
\begin{cases}
  e_0= y+\pi\bar x y
  &\\
  e_{\alpha}= \pi \bar x y+ \sum_{i=0}^{p^{n_1}-1} \alpha^{i} x^{(i)}y &
  \text{ for } \alpha\in\F_{p^{n_1}}^{\ast}\cdot\sigma,
\end{cases}
\end{equation*}
where $e_\alpha$ corresponds to the eigenvector $\alpha$.
The former formula can be seen as a special case of the latter if we agree to set $0^0=1$.
See~\cite[Section~2]{CaMa:Hamiltonian} for some comments on how to obtain and conveniently manipulate formulas like these.

The decomposition of $L$ into weight spaces
with respect to $\ad e_0$ extends the grading of $W(1;n_1)$ to a grading of $L$.
Since $e_0$ belongs to $\bar L_0$ it normalizes every component $\bar L_k$.
Intersecting this grading of $L$ with the $\Z$-grading $\bigoplus_{k\in\Z}\bar L_k$
yields a $(\Z \times \F_{p^{n_1}})$-grading, from which we will obtain the desired cyclic grading of $L$.
Thus, our first goal is finding the weight spaces of $\ad e_0$ on each component $L_j$.
We have
\begin{align*}
  \{e_0, x^{(i)} y^{(j)}\} &= \{y+\pi\bar x y, x^{(i)}
  y^{(j)}\}=\\
  & =
\begin{cases}
  (1-j)\bar xy^{(j)}
  - j\pi x^{(p^{n_1}-2)} y^{(j)}
  &\text{if $i=0$,}\\
  y^{(j)}+ (1+j)\pi \bar x y^{(j)}
  &\text{if $i=1$,}\\
  x^{(i-1)} y^{(j)} &\text{if $i>1$.}
\end{cases}
\end{align*}

Choose $\rho\in\F$ such that $\rho^{p^{n_1}}-\pi\rho-1=0$.
It is not hard to see that the characteristic polynomial for $\ad e_0$ on $\bar L_{1-j}$ is
$x^{p^n_1}-\pi x+j-1$, and its roots form the set
$(1-j)\rho+\F_{p^{n_1}}\cdot\sigma$.
Eigenvectors for $\ad e_0$ on $\bar L_{1-j}$, for
$j=0,\ldots,p^{n_2}-1$, are
\begin{equation}\label{eq:e_jalpha}
  e_{1-j,\alpha}=
  \pi j\bar x y^{(j)}+ \sum_{i=0}^{p^{n_1}-1}
  \alpha^{i}
  x^{(i)}y^{(j)}
  \qquad \text{for $\alpha\in(1-j)\rho+\F_{p^{n_1}}\cdot\sigma$}.
\end{equation}
These eigenvectors extend the complete set of eigenvectors for $\ad e_0$ on the subalgebra $W(1;n_1)$
given by its basis elements $e_{0,\alpha}=e_{\alpha}$.
We have chosen our notation so that
$e_{1-j,\alpha}\in L_{1-j}$, and
\begin{equation*}
\{e_0,e_{1-j,\alpha}\}=
\alpha e_{1-j,\alpha}.
\end{equation*}
Note that, differently from~\cite[Section~7]{CaMa:Hamiltonian}, the subscript $1-j$ is a genuine integer
and is not to be read modulo $p^{n_2}-1$.
Also, the two subscripts are not independent, since the value of $1-j$ modulo $p$ is determined by $\alpha$.
Thus, we have obtained a grading
$L=\bigoplus_{(k,\alpha)} L_{k,\alpha}$
over a certain subgroup of $\Z\times(\F_p\cdot\rho+\F_{p^{n_1}}\cdot\sigma)$
isomorphic with
$\Z\times\F_{p^{n_1}}$,
the degree of each basis element $e_{k,\alpha}$ being given by its pair of subscripts $(k,\alpha)$.
Note that all components $L_{k,\alpha}$ of the grading have dimension zero or one,
the latter possibility occurring for $2-p^{n_2}\le k\le 1$.

It is easy to find the complete multiplication table of the basis $e_{1-j,\alpha}$.
We have
\begin{equation}\label{eq:multiplication}
\{e_{1-j,\alpha},e_{1-l,\beta}\}=
\left(\beta\binom{j+l-1}{l}-\alpha\binom{j+l-1}{j}\right)
e_{2-j-l,\alpha+\beta},
\end{equation}
which is to be read as zero when $2-j-l$ lies outside its legal range
$2-p^{n_2},\ldots,0,1$.
In fact, since the basis is graded we know that the result on the right-hand side
must be a scalar multiple of $e_{2-j-l,\alpha+\beta}$.
Thus, it remains only to determine that scalar,
for example by computing the coefficient of the monomial $y^{(j+l-1)}$ in the result
(noting that $y^{(j)}$ always appears with nonzero coefficient in $e_{1-j,\alpha}$).
To accomplish this it suffices to compute the product of the only relevant terms, namely
\[
\{y^{(j)},\beta xy^{(l)}\}+\{\alpha xy^{(j)},y^{(l)}\}=
\left(\beta\binom{j+l-1}{l}-\alpha\binom{j+l-1}{j}\right)
y^{(j+l-1)}.
\]

Now we restrict our attention to the case of main interest for us by setting $n_1=1$.
For later reference we rewrite Equation~\eqref{eq:e_jalpha} in this special case:
\begin{equation}\label{eq:e_jalpha-small}
  e_{1-j,\alpha}=
  \sigma^{p-1} j\bar x y^{(j)}+ \sum_{i=0}^{p-1}
  \alpha^{i}
  x^{(i)}y^{(j)}
  \qquad \text{for $\alpha\in(1-j)\rho+\F_{p}\cdot\sigma$}.
\end{equation}
We introduce the shorthand $q=p^{n_2}$, hence $L$ has dimension $pq$.
Roughly speaking, we want to produce a cyclic grading of $L$
by ``viewing $j$ modulo $q-1$'', but we have to be careful because $j$ and $\alpha$ are not independent.
In precise terms this means factoring the grading group, which is the subgroup
$\langle (1,\rho),(0,\sigma)\rangle$ of
$\Z\times(\F_p\cdot\rho+\F_{p}\cdot\sigma)$, modulo its subgroup
$\langle(1-q,\rho)\rangle$.
We naturally obtain a new grading of $L$ over the quotient group $G$, which is cyclic of order $(q-1)p$.
We need not introduce a specific notation for this grading, but for the loop algebra construction
to make sense we must choose an isomorphism of $G$ with $\Z/(q-1)p\Z$.
We choose the isomorphism which maps
$(1,\rho+\sigma)+\Z\cdot(1-q,\rho)$ to $1+(q-1)p\Z$.
In other words, we choose
the former as the distinguished generator ``$\bar 1$'' of $G$
referred to in Section~\ref{sec:Nottingham}.
The corresponding component of the cyclic grading is
$L_{1,\rho}\oplus L_{2-q,2\rho+\sigma}$.

Next we prove that the corresponding loop algebra is thin with second diamond in degree $q$,
thus a Nottingham Lie algebra,
and that all diamonds are of finite type.
All components in the $(\Z/(q-1)p\Z)$-grading of $L$ are one-dimensional, except those
of degree congruent to $1$ modulo $q-1$, which are two-dimensional
(having the form $L_{1,\alpha}+L_{2-q,\alpha+\rho}$)
and will give rise to the diamonds in the loop algebra.
We proceed as in Section~\ref{sec:mixed}
and set $X=e_{1,\rho+\sigma}$ and $Y=e_{2-q,2\rho+\sigma}$.
Equation~\eqref{eq:multiplication} yields the formulas
\begin{align*}
  & \{e_{1-j,\alpha},X\}=(\rho+\sigma) e_{2-j,\alpha+\rho+\sigma}
  \quad\text{for $j=0,\ldots, q-1$},\\
  & \{e_{1-j,\alpha},Y\}=0 \quad\text{for $j=2,\ldots q-1$},\\
  & \{e_{0,\alpha},Y\}=
  (\alpha+2\rho+\sigma) e_{2-q,\alpha+2\rho+\sigma},\\
  & \{e_{1,\alpha},Y\}=
  -\alpha e_{3-q,\alpha+2\rho+\sigma},
\end{align*}
which describe the adjoint action of $X$ and $Y$ on $L$.
Note that all coefficients which appear are never zero.
Hence, the formulas show at once that
\begin{align*}
&\{L_{k,\alpha},\F X+\F Y\}=
L_{k+1,\alpha+\rho}
\quad\textrm{for $k\neq 2-q,0,1$, and}
\\
&\{L_{0,\alpha},\F X+\F Y\}=
L_{1,\alpha+\rho}+L_{2-q,\alpha+2\rho},
\end{align*}
which imply that the covering property holds in these degrees, since $L_{k,\alpha}$ is one-dimensional.
To see that the loop algebra is thin, it remains to check the covering
property on the two-dimensional components.
The elements $\{e_{0,\alpha},X\}$ and $\{e_{0,\alpha},Y\}$ span the component
$L_{1,\alpha+\rho+\sigma}+L_{2-q,\alpha+2\rho+\sigma}$,
and we have
\begin{align*}
  & \{e_{0,\alpha},X,X\}= 0,\\
  & \{e_{0,\alpha},X,Y\}=
  -(\rho+\sigma)(\alpha+\rho+\sigma) e_{3-q,\alpha+3\rho+2\sigma},\\
  & \{e_{0,\alpha},Y,X\}=
  (\rho+\sigma)(\alpha+2\rho+\sigma)e_{3-q,\alpha+3\rho+2\sigma},\\
  & \{e_{0,\alpha},Y,Y\}= 0.
\end{align*}
The first and fourth equations together imply that the covering property holds.
Thus, the loop algebra of $L$ is a thin Lie algebra, with second diamond in degree $q$.
The second and third equations determine the diamond types.
In fact, the element $V_{t}= e_{0,-(t-1)\sigma}$ (which we may define for every positive integer $t$)
gives rise to a nonzero element in the component preceding the $t$th diamond in the loop algebra,
and satisfies
\begin{equation*}
(1-\mu_t)\{V_{t},X,Y\}=
\mu_t \{V_{t},Y,X\},
\end{equation*}
with
$\mu_t=-1+(t-2)\sigma/\rho$.
Consequently, the $t$th diamond in the loop algebra,
which occurs in degree $(t-1)(q-1)+1$
and originates from the component of $L$ spanned by
$\{V_{t},X\}$ and $\{V_{t},Y\}$,
has type $\mu_t$.
We have completed a proof of the following result.

\begin{theorem}\label{thm:grading_finite}
Let $\F$ have odd characteristic $p$, and let $q=p^n$, where $n>0$.
Assume that there are $\sigma,\rho\in\F$ with $\sigma\neq 0$ and $\rho^p-\sigma^{p-1}\rho-1=0$.
A grading of $H(2;(1,n);\Phi(1))$ over the integers modulo $(q-1)p$ is determined by declaring
$X=e_{1,\rho+\sigma}$ and $Y=e_{2-q,2\rho+\sigma}$, that is,
\[
X=\sum_{i=0}^{p-1}
(\rho+\sigma)^{i}x^{(i)}
\qquad\textrm{and}\qquad
Y=-\sigma^{p-1}\bar x\bar y+
\sum_{i=0}^{p-1}
(2\rho+\sigma)^{i}x^{(i)}\bar y,
\]
to have degree one.
They generate $H(2;(1,n);\Phi(1))$ and are part of a graded basis
given by the elements in~\eqref{eq:e_jalpha-small}, where $\pi=\sigma^{p-1}$.
The degree $k$ modulo $(q-1)p$ of $e_{r,r\rho+s\sigma}$ is determined by
$k\equiv r\pmod{q-1}$ and $k\equiv s\pmod{p}$.
The corresponding loop algebra $N$ is thin with second diamond in degree $q$.
The diamonds occur in all degrees congruent to $1$ modulo $q-1$,
and their types follow an arithmetic progression not contained in the prime field.
\end{theorem}

Since the third diamond has always type $-1$, the arithmetic progression which gives the diamond types
is determined by the types of the third diamond, which is $\mu_3=-1+\sigma/\rho$.
Conversely, if the type $\mu_3$ of the third diamond is assigned, with $\mu_3\in\F\setminus\F_p$,
then $\sigma$ and $\rho$
are uniquely determined by
\begin{equation}\label{eq:assigned-type}
\sigma^p\left(\frac{1}{\mu_3^p+1}-\frac{1}{\mu_3+1}\right)=1
\qquad\text{and}\qquad
\rho(\mu_3+1)=\sigma.
\end{equation}
Consequently, Theorem~\ref{thm:grading_finite} gives an existence proof for the Nottingham Lie algebras
whose uniqueness was established in~\cite[Theorem~2.3]{CaMa:Nottingham} for $p>5$.
In essence, that result stated that up to isomorphism there is a unique (infinite-dimensional) thin Lie algebra
whose largest quotient of class $2q$ is isomorphic with $L/L^{2q+1}$,
where $N$ is the thin Lie algebra constructed in Theorem~\ref{thm:grading_finite}.
An equivalent formulation is that $N$ is the unique Nottingham Lie algebra with second diamond in degree $q$,
and third diamond in degree $2q-1$ and prescribed type $\mu_3\in\F\setminus\F_p$ (for $p>5$).
Note that distinct values for $\mu_3$ correspond to non-isomorphic Nottingham Lie algebras.
In fact, the subspaces $\langle X\rangle$ and $\langle Y\rangle$ of $L_1$ are invariants of the algebra,
as noted after the proof of~\eqref{eq:second_diamond_type}, and $\mu_3$ is unaffected by replacing $X$ and $Y$
with scalar multiples.

Equations~\eqref{eq:assigned-type} also justify a claim which we made at the beginning of this section.
In fact, they show that the ground field $\F$ satisfies the assumptions of
Theorem~\ref{thm:grading_finite} provided it is perfect and contains what we want to be the type
of the third diamond of the resulting loop algebra $N$,
namely, $\mu_3=-1+\sigma/\rho$.
However, note that the existence of a thin algebra with structure as described in~\cite[Theorem~2.3]{CaMa:Nottingham}
does not depend on the perfectness assumption on $\F$.
This is because~\cite[Theorem~2.3]{CaMa:Nottingham} and its proof show that the loop algebra
$N$ of Theorem~\ref{thm:grading_finite} is actually defined over the field $\F_p(\mu_3)$.
Thus, when $\F$ is not perfect $N$ may be a loop algebra of an $\F$-form of $H(2;(1,n);\Phi(1))$
different from the one considered here.

\begin{rem}\label{rem:char_two_finite}
We can interpret Theorem~\ref{thm:grading_finite} in characteristic two in a similar way as we did for
Theorem~\ref{thm:grading_mixed} in Remark~\ref{rem:char_two_mixed}.
Note first that in odd characteristic the only basis elements in~\eqref{eq:e_jalpha-small} which involve $\bar x\bar y$ are
$e_{2-q,\alpha}$ for $\alpha\in 2\rho+\F_{p}\cdot\sigma$.
However, in characteristic two we have $e_{2-q,\sigma}=\bar y$, and so the only basis element involving $\bar x\bar y$
is $e_{2-q,0}=\sigma\bar x\bar y+\bar y$.
Thus, all the remaining basis element span the derived subalgebra
$H(2;\n;\Phi(1))^{(1)}$.
In fact, the elements $X=1+(\rho+\sigma)x$ and $Y=\bar y$ generate
$H(2;\n;\Phi(1))^{(1)}$.
The loop algebra $N$ has a fake diamond of type $-1\equiv 1$ in each degree congruent to $q$ mod $2q-2$,
and a genuine diamond of type $\mu_3=-1+\sigma/\rho$ in each degree congruent to $1$ mod $2q-2$.
Together with the Nottingham algebra described in Remark~\ref{rem:char_two_mixed},
this Nottingham algebra $N$ falls into one of the cases described in~\cite{JuYo:quotients},
where the largest quotient of maximal class of $N$ over a term of the lower central series
is not metabelian, but has constituent sequence $2q',2q'-2$.

In the special case where $p=q=r=2$ the algebra
$H(2;\n;\Phi(1))^{(1)}$ is the three-dimensional simple Lie algebra
and the construction of Theorem~\ref{thm:grading_finite}
yields yet another realization of the thin Lie algebra considered in~Remark~\ref{rem:char_two_mixed_special}.
\end{rem}

\section{Nottingham Lie algebras with diamonds of finite types in the prime field}\label{sec:prime_field}

Theorem~\ref{thm:grading_finite} produces a Nottingham Lie algebra with second diamond in degree $q$
and third diamond in degree $2q-1$ and of type $\mu_3$, for any prescribed $\mu_3\not\in\F_p$.
As we have mentioned in the Introduction,
analogous Nottingham Lie algebras with $\mu_3\in\F_p$ were constructed in~\cite{Car:Nottingham,Avi}
as loop algebras of $W(n;1)$ if $\mu_3=-1$, and $H(2;(1,n))$ otherwise (suitably extended by an outer derivation
if $\mu_3=-2$ or $-3$).
In this section we show how the construction
of Section~\ref{sec:finite} allows us to recover those Nottingham Lie algebras as well.

We deal first with the easier case $\mu_3=-1$.
In the setting of Section~\ref{sec:finite} we may interpret this as $\sigma=0$, whence $\rho^p-1=0$, and hence $\rho=1$.
The element $e_0$ defined as in Section~\ref{sec:finite} is now $e_0=y$.
Differently from the original setting of Section~\ref{sec:finite}, here $e_0$ is not semisimple.
In fact, the only admissible value for $\alpha$ in Equation~\eqref{eq:e_jalpha-small} is $\alpha=(1-j)\rho$,
and the formula yields the only eigenvector $e_{1-j,\alpha}$ for $\ad e_0$ on $\bar L_{1-j}$.
The collection of these eigenvectors
spans a Lie subalgebra of
$L$ of dimension $q$.
In fact, the multiplication formula given in Equation~\eqref{eq:multiplication}
still applies in this degenerate situation and yields, after a simplification,
\[
\{e_{1-j,(1-j)\rho},e_{1-l,(1-l)\rho}\}=
\rho\left(\binom{j+l-1}{l}-\binom{j+l-1}{j}\right)
e_{2-j-l,(2-j-l)\rho}.
\]
It follows that the elements $X$ and $Y$ as defined in Theorem~\ref{thm:grading_finite}
generate a subalgebra of $H(2;(1,n);\Phi(1))$ isomorphic with $W(1;n)$
via $E_i\mapsto -\rho^{-1}e_{-i,-i}$.
The grading of $W(1;n)$ thus defined, over the integers modulo $q-1$, coincides
with that used in~\cite{Car:Nottingham}
to construct a Nottingham Lie algebra with all diamonds of type $-1$.
(A similar situation would occur in the contruction given in~\cite[Section~5]{Avi}:
setting $\mu_3=-1$ in that context, a value excluded there, then
$X$ and $Y$ would generate a subalgebra of $H(2;(1,n))$ isomorphic with $W(1;n)$.)

In the remaining cases $\mu_3\in\F_p\setminus\{-1\}$, the ratio $\sigma/\rho$ should also belong to the prime field,
in contradiction with the assumption $\rho^p-\sigma^{p-1}\rho-1=0$ of Theorem~\ref{thm:grading_mixed}.
However, we can make sense of these cases in terms of deformations.
Loosely speaking, we can insert a parameter $\varepsilon$ in the Lie bracket defined by
Equations~\eqref{eq:Poisson_1} and~\eqref{eq:Poisson_exception} leaving the isomorphism type
of the algebra unchanged as long as $\varepsilon\neq 0$.
Then Theorem~\ref{thm:grading_finite} gives a grading of the algebra where the diamond types of the corresponding
loop algebra depend on $\varepsilon$.
It turns out that when $\varepsilon$ ``approaches zero'' the type of the third diamond
``approaches an element of the prime field''.
However, the isomorphism type of the finite dimensional algebra also changes in the limit.
It is possible to make this rigorous by adopting the language of Gerstenhaber's
deformation theory of algebras~\cite{Ger:deformation1,GerSch}.
However, to avoid the risk of making this section less readable and distracting from the essence of the argument
we have preferred an informal exposition,
thus relegating the technicalities to Remark~\ref{rem:deformation}.

Assume for simplicity that $\F$ is algebraically closed.
Let $\varepsilon$ be an element of $\F$
and define a multiplication $\{\ ,\ \}_\varepsilon$ on $\mathcal{O}(2;(1,n))$
by the same formula in Equation~\eqref{eq:Poisson_1} for $i+k>0$, and
\begin{equation}\label{eq:Poisson_epsilon}
\{y^{(j)},y^{(l)}\}_\varepsilon=
\varepsilon\left(
\binom{j+l-1}{l}-
\binom{j+l-1}{j}
\right)
\bar xy^{(j+l-1)}
\end{equation}
replacing Equation~\eqref{eq:Poisson_exception}.
This amounts to replacing the automorphism $\Phi(1)$ of $\mathcal{O}((2))$ with the automorphism
$\Phi(1)_\varepsilon$ such that
$\Phi(1)_\varepsilon(x)=x\exp(\varepsilon x^{(p^{n_1})})$
and $\Phi(1)_\varepsilon(y)=y$.
The resulting Lie algebra $H(2;(1,n);\Phi(1)_\varepsilon)$
is isomorphic with $H(2;(1,n);\Phi(1))$
if $\varepsilon\neq 0$, and with a central extension of $H(2;(1,n))^{(1)}$ if $\varepsilon=0$,
which we naturally identify with $\mathcal{O}(2;(1,n))$ with the Poisson bracket~\eqref{eq:Poisson_0},
and temporarily call $\widehat H$ for convenience.
Let $\sigma$ be a nonzero element of $\F$, and set $e_0=y+\sigma^{p-1}\bar xy$.
The characteristic polynomial for $\ad e_0$ on $\bar L_{1-j}$ is
$x^{p}-\sigma^{p-1} x+\varepsilon(j-1)$.
Its roots consist of the set
$(1-j)\rho+\F_p\cdot\sigma$, where we have chosen $\rho$ so that $\rho^p-\sigma^{p-1}\rho-\varepsilon=0$.
Eigenvectors for $\ad e_0$ on $\bar L_{1-j}$,
for $j=0,\ldots,q-1$, are still given by Equation~\eqref{eq:e_jalpha-small}
and multiply according to Equation~\eqref{eq:multiplication}.

When we set $\varepsilon=0$, Equation~\eqref{eq:e_jalpha-small} gives a complete set of eigenvectors for
$\ad e_0$ in $\widehat H$, where $e_0=y+\sigma^{p-1}\bar xy$ (as long as $\sigma\neq 0$).
Now $\rho$ and $\sigma$ are related by the equation $\rho^p-\sigma^{p-1}\rho=0$.
When $\rho=0$, which should correspond to the third diamond having type $\mu_3=\infty$,
this grading does not yield a thin algebra
(because $\{e_{1,0},Y\}=0=\{e_{0,-\sigma},Y\}$).
However, we have already produced various Nottingham algebras with third diamond of type $\infty$ in Section~\ref{sec:mixed}
as loop algebras of $H(2;\n;\Phi(1))$,
and hence we assume $\rho\neq 0$.
Since the diamond types depend only on the ratio $\sigma/\rho$, which is now a nonzero element of the prime field,
we may as well assume that $\sigma$ and $\rho$ belong to the prime field, whence $\sigma^{p-1}=1$
produces a slight simplification in formula~\ref{eq:e_jalpha-small}.
Assigning degrees to the elements $e_{1-j,\alpha}$ with the same rule used in Theorem~\ref{thm:grading_finite}
gives a grading of $\widehat H$ over the integers modulo $p(q-1)$.
Before we consider the corresponding loop algebra,
recall that $\widehat H$ is not simple, unlike $H(2;(1,n);\Phi(1))$.
Its subalgebra
$\langle x^{(i)}y^{(j)}:i+j<p+q-2\rangle$
of codimension one is spanned by the collection of all $e_{1-j,\alpha}$ with the exception of
$e_{2-q,0}=\bar y-\bar x\bar y$.
Furthermore, $\widehat H$ has a one-dimensional centre spanned by $e_{1,0}=1$, and the quotient
$\langle x^{(i)}y^{(j)}:i+j<p+q-2\rangle/\langle 1\rangle$ is a simple Lie algebra
(namely, $H(2;(1,n))^{(2)}$).

It is easy to verify that the elements $X=e_{1,\rho+\sigma}$ and $Y=e_{2-q,2\rho+\sigma}$
generate $\widehat H$ provided $\sigma/\rho\neq -1$.
The corresponding loop algebra $L$ is not thin because it has a nonzero centre.
However, the quotient of $L$ modulo its centre is a Nottingham Lie algebra
and has diamonds, possibly {\em fake}, in all degrees congruent to $1$ modulo $q-1$.
Their types follow an arithmetic progression contained in the prime field, determined by the type
$-1+\sigma/\rho$ of the third diamond.
Up to minor differences in notation, this matches the construction of these algebras given in Theorem~3 of~\cite{Avi}.

In the case $\sigma/\rho=-1$ which was excluded above, $X$ turns out to be the central element
$e_{1,0}$ of $\widehat H$.
In fact, the arithmetic progression of diamond types would predict that the very first diamond in the corresponding
loop algebra is fake of type $0$.
This is one of two cases in~\cite{Avi} where to produce a thin algebra one needs to introduce an outer derivation
of the simple algebra, in this case $\ad x^{(p)}$.
The other exceptional case in~\cite{Avi} corresponds to our case $\sigma/\rho=-2$, where the arithmetic progression
of diamond types would predict the first diamond being fake of type $1$.
There one needs to introduce the outer derivation $\ad\bar x\bar y$ of the simple Hamiltonian algebra,
which was already taken care of by the present construction.

\begin{rem}\label{rem:deformation}
As promised earlier in this section, we briefly sketch how one can make the above deformation argument rigorous.
If we let $\varepsilon$ be an indeterminate over $\F$ and denote by $\F(\varepsilon)$
the corresponding field of rational functions
then Equations~\eqref{eq:Poisson_1} and~\eqref{eq:Poisson_epsilon} define a multiplication
$\{\ ,\ \}_\varepsilon$ on the divided power algebra $\mathcal{O}(2;(1,n))$ over $\F(\varepsilon)$.
Because this Lie algebra is actually defined over the polynomial ring $\F[\varepsilon]$,
one obtains a Lie algebra over $\F$ by specializing $\varepsilon$ to any element of $\F$.
However, in order to proceed further and deal with the eigenvectors for $\ad e_0$
it is necessary to enlarge $\F(\varepsilon)$ to a field containing $\sigma$ and $\rho$.
The easiest way to do this is to introduce a further indeterminate $\sigma$
and view the Lie algebra over the field $\E=\F(\sigma)((\varepsilon))$ of formal Laurent series over $\F(\sigma)$.
The Lie algebra then is (a form of) the Hamiltonian algebra $H(2;(1,n);\Phi(1))$ over $\E$.
The polynomial $Z^p-\sigma^{p-1}Z-\varepsilon$ has all its roots $\rho_0,\rho_0+\sigma,\ldots,\rho_0+(p-1)\sigma$ in $\E$,
where $\rho_0=-\sigma\sum_{i= 0}^{\infty}(\varepsilon/\sigma^p)^{p^i}$.
Taking any of the roots as $\rho$, Equation~\eqref{eq:e_jalpha-small} gives a complete set of eigenvectors for $\ad e_0$.
The loop algebra constructed according to the recipe given in Theorem~\ref{thm:grading_finite}
is a Nottingham algebra, over the field $\E$, with all diamonds of finite type not all in the prime field.
Since everything seen so far is actually defined over the ring $\F(\sigma)[[\varepsilon]]$ of formal power series,
we are allowed to specialize $\varepsilon$ at zero.
The Lie algebra then becomes a central extension $\widehat H$ of
$H(2;(1,n))^{(1)}$ over $\F(\sigma)$,
the element $\rho$ becomes one of the roots $0,\sigma,\ldots,(p-1)\sigma$ of
$Z^p-\sigma^{p-1}Z$, and the corresponding loop algebra has all diamond in the prime field, as described above.
In Gerstenhaber's language of deformations~\cite{Ger:deformation1,GerSch} one says that the Lie algebras (both the Hamiltonian algebra
and its loop algebra) over $\F(\sigma)[[\varepsilon]]$ are {\em formal deformations}
of the corresponding ones over $\F(\sigma)$.
\end{rem}

\bibliography{References}

\end{document}